\documentclass[11pt,letterpaper]{article}
\usepackage{times}
\usepackage{amsmath,amsfonts,amssymb}

\newtheorem{theorem}{Theorem}
\newtheorem{lemma}{Lemma}[section]

\newcommand{\remark}[1][Remark]{\vspace{1ex} \noindent \textbf{#1: }}

\newcommand{\Z}{\mathbb{Z}}
\newcommand{\R}{\mathbb{R}}
\newcommand{\T}{\mathbb{T}}
\newcommand{\Rd}{\R^d}
\newcommand{\RR}{\mathbb{R}}
\newcommand{\ip}[2]{{\langle#1,#2\rangle}}
\newcommand{\cat}{{:}}
\newcommand{\N}{{\mathbb N}}
\newcommand{\B}{\mathcal{B}}
\newcommand{\A}{\mathcal{A}}

\newcommand{\E}{\mathsf E}
\newcommand{\PP}{\mathbb{P}}

\newcommand{\eps}{\varepsilon}
\newcommand{\sgm}{\sigma}
\newcommand{\ONE}{\mathbf{1}}
\newcommand{\ccdot}{ \ \cdot \ }
\newcommand{\Norm}[1]{\lvert\!\lvert\!\lvert #1 \rvert\!\rvert\!\rvert}
\newcommand{\Nice}[1][\rho]{\mathcal{N}_{#1}}
\newcommand{\bpf}[1][Proof]{{\noindent {\sc #1: }}}
\newcommand{\epf}{{{\hspace{4 ex} $\Box$ \smallskip}}}

\newcommand{\Pl}{\Pi_\ell}
\newcommand{\Ph}{\Pi_h}
\newcommand{\DV}{\mathcal{D}(V)}
\newcommand{\Dc}{\mathcal{D}}
\newcommand{\DU}{\mathcal{D}(U)}
\newcommand{\LL}{\mathbb{L}}
\newcommand{\HH}{\mathbb{H}}
\newcommand{\X}{\mathbb{X}}

\title{Stationary Solutions of Stochastic Differential Equation with
  Memory and Stochastic Partial Differential Equations.}
\author{Yuri Bakhtin\thanks{School of Math, Institute for Advanced Study,
    Princeton, NJ 08540, USA and International Institute for Earthquake
Prediction Theory and Mathematical Geophysics, Warshavskoye sh. 79, kor. 2,
Moscow 113556,  Russia
 }
\and Jonathan C. Mattingly
\thanks{School of Math, Institute for Advanced Study,
    Princeton, NJ 08540, USA and Department of Math, Duke University,
  Durham, NC 27708, USA} }
\date{July 2nd, 2003}
\begin{document}
\maketitle
\begin{abstract}
  We explore It\^o stochastic differential equations where the drift
  term possibly depends on the infinite past. Assuming the
  existence of a Lyapunov function, we prove the existence of a stationary
  solution assuming only minimal continuity of the
  coefficients. Uniqueness of the stationary solution is proven if
  the dependence on the past decays sufficiently fast. The results of
  this paper are then applied to stochastically forced dissipative
  partial differential equations such as the stochastic Navier-Stokes
  equation and stochastic Ginsburg-Landau equation.

  \vspace{1ex}
\end{abstract}

\section{Introduction}

This note explores the ergodic theory of It\^o stochastic differential
equations with memory. Specifically, we consider equations on $\RR^d$
with additive noise of the form
\begin{equation}
dX(t) = a(\pi_tX)dt + dW(t).
\label{sdde}
\end{equation}
Here $W(t),t\in\R$ is a standard $d-$dimensional Wiener process
(i.e. a Gaussian $\Rd$-valued stochastic process
with continuous trajectories defined on the whole real line $\R$ with
independent and stationary increments, $W(0)=0$, $\E W(t) = 0$,
and $\E W_i(t)W_j(t)=|t|\delta_{ij}, t\in\R$, $i,j=1,\ldots,d$).
The projection and shift $\pi_t$ is a map from the space
$C$ of $\Rd$-valued continuous functions defined on $\R$ to the space
$C^-$ of continuous functions defined on $\R_-=(-\infty,0]$. $\pi_t$
is defined by
\begin{equation*}
  \pi_tX(s)=X(s+t),\quad s\in\R_-.
\end{equation*}
This map gives the past of a continuous process up to time
$t\in\R$. Lastly the map $a:C^- \rightarrow \RR^d$ gives the effect of
the past on the present moment of time.

The first results on stationary solutions for SDEs of this type
appeared in \cite{b:ItoNisio}.  In this paper, we give new sufficient
conditions for existence and uniqueness of stationary solution for
equation (\ref{sdde}) which may be considered as non-Markovian
(``Gibbsian'') counterparts of those in \cite{b:Veretennikov97}. Here
we do not address the question of mixing. However \cite{b:Mattingly02}
provides the needed framework to address this question in the
non-Markovian setting.

We hope that our results will be useful in many different contexts.
However, one of the main guiding examples has been recent progress in
the ergodic theory of stochastic partial differential equations driven
by white noise found in
\cite{b:BricmontKupiainenLefevere01,b:EMattinglySinai00,b:ELui02,b:Mattingly02}.
In \cite{b:KuksinShirikyan00} and subsequent papers related ideas were
developed independently in kicked noise setting. In the
\cite{b:MasmoudiYoung02} the kicked setting is furtherer developed. In
particular, we follow the ideas as laid out in
\cite{b:EMattinglySinai00}.  There using ideas of determining modes
and inertial manifolds (see \cite{b:FoiasProdi67, b:Te88}), the
infinite dimensional diffusion is reduced to an It\^o process with
memory on a finite dimensional phase space. If the resulting It\^o
process is elliptic then the methods of this note can be used to
establish uniqueness of stationary measure.  In
\cite{b:Mattingly02,b:Hairer02}, similar ideas were used to control a
convergence rate to the invariant measure; there the packaging was
more Markovian. These ideas were further described and extended in
\cite{b:KuksinShirikyan02} and \cite{b:Mattingly03Pre}.  With the
exception of \cite{b:EMattinglySinai00,b:Bakhtin02} the memory is less
explicit in the preceding works, here we bring it to the foreground
and give general conditions under which the ideas are applicable. We
emphasize that we do not prove fundamentally new ergodic results for
stochastic PDEs. The results in
\cite{b:Mattingly02,b:Hairer02,b:KuksinShirikyan02,b:Mattingly03Pre}
cover our PDE examples and give convergence rates in addition. However, we
wish to clarify the reduction to an SDE with memory as it provides
useful intuition. We believe that our results are new in the context
of general SDEs with memory. Lastly similar ideas to those in this
note can be used to treat equations with state dependent diffusion
matrices, provided that the system is uniformly elliptic.

After developing the general theory, we examine a number of
pedagogical examples. Then we show how one can use the results of this
paper to reduce a dynamical system to one of smaller dimension but
with memory by removing stable dimensions. This reduced system can
then be analyzed to understand the asymptotic properties of the
system.  The construction reflects the simple fact that the long term
dynamics is dictated by the dynamics in the unstable directions.  This
reduction is particularly useful in the context of dissipative partial
differential equations where the reduced system is finite dimensional;
and hence, it has a simpler topological structure than the original
equation. In that setting, it is our imperfect understanding of the fine scale
topological structure in the Markovian setting which limits our
progress. By switching to the finite dimensional setting where
all relevant topologies are equivalent and the Lebesgue measure exists, we
can make progress.

\section{Definitions and Main Results}\label{sec:Definitions}
Along with the set of pasts $C^-$ defined in the introduction, we also
define the set of futures $C^+$ as the space of $\RR^d$-valued
continuous functions on $\R_+=[0,\infty)$. We denote by $\pi_+$ the
natural projection from $C \rightarrow C^+$. For a stochastic process
$X$ and a set $A\subset\R$ the $\sgm$-algebra generated by random variables
$X(s), s\in A$ will be denoted by $\sgm_A(X)$ and the $\sgm$-algebra
generated by random variables $X(s)-X(t), s,t\in A$ will be denoted by
$\sgm_A(dX)$.
Consider the space $\Omega=C\times C$ with LU-topology defined by the metric
\begin{equation*}
  \rho(f,g)=\sum_{n=-\infty}^\infty 2^{-|n|}(\|f-g\|_n\wedge1),\quad
  f,g\in C\times C,
\end{equation*}
where $\|(h_1,h_2)\|_n=\max_{-n\le t\le n}(|h_1(t)|+|h_2(t)|)$ for
$(h_1,h_2)\in \Omega$ and $|\cdot|$ denotes the Euclidean norm.

\remark[Solutions of the SDE]A probability measure $P$ on the space
$\Omega$ with Borel $\sgm$-algebra $\B$ is said to define a
\textsc{solution} to the equation \eqref{sdde} on some subset
$\mathcal{R}$ of $\R$ if the following three conditions are fulfilled
with respect to the measure $P$:
\begin{enumerate}
\item The projection $W:C\times C\to C$, $\omega = (\omega_1,\omega_2)\mapsto\omega_2$,
is a standard $d$-dimensional Wiener process.
\item Henceforth, let $X$ denote the project $X:C\times C\to C$,
  $\omega = (\omega_1,\omega_2)\mapsto\omega_1$. For any $t\in\R$
\begin{equation}
  \sgm_{(-\infty,t]}(X)\vee\sgm_{(-\infty,t]}(dW)\quad \mbox{is independent of}\quad
  \sgm_{[t,\infty)}(dW).
  \label{independence}
\end{equation}
\item If $s<t$ and $s,t\in\mathcal{R}$ then
  \begin{equation}
    X(t)-X(s)\stackrel{\mbox{\small a.s.}}{=}\int_s^ta(\pi_\theta X)d\theta+W(t)-W(s).
    \label{fulfilsdde}
  \end{equation}
\end{enumerate}
The process $X$ or the couple $(X,W)$ will be also often
referred to as solution for equation \eqref{sdde}.

\remark[Stationary Solutions]If $P$ defines a solution on $\R$ and
in addition the distribution of the process
\begin{equation*}
  (X,dW) \equiv (X(t), -\infty<t<\infty, W(v)-W(u), -\infty<u<v<\infty)
\end{equation*}
does not change under time shifts then the measure $P$ is said to define a
\textsc{stationary solution}.

\remark[Solution to Cauchy Problem] We will always assume strong existence and
pathwise uniqueness of solution to the Cauchy problem with initial data which grows
sufficiently slow at $-\infty$.

More precisely, for $\rho>0$ denote $C^-_\rho$ the space of trajectories $x\in C^-$ for which
\begin{equation*}
\|x\|_\rho=\sup_{t\in\R_-}\frac{|x(t)|}{1+|t|^\rho}<\infty.
\end{equation*}
We will require that for some $\rho>0$ and any $x\in C^-_\rho$, which possesses
certain averaging property (which will be described after the definition of
Lyapunov function below)
there exists
a measurable map $\Phi:C^+\rightarrow C^+$
such that if $W$ is a standard Wiener process
under some measure $P$
then $(X,W)$ is a solution to \eqref{sdde} on $\R_+$
where
$X(t)=x(t)\ONE_{\{t\le0\}}+\Phi(\pi_+ W)(t)\ONE_{\{t>0\}}$
which means that relations \eqref{independence}
and \eqref{fulfilsdde} are true on $\R_+$.
Moreover,
if under some measure $P$ the process $(X,W)$ is a solution to
\eqref{sdde} on $\R_+$ and $X(t)=x(t)$ for $t\le0$ then
$X(t)=\Phi(\pi_+ W)(t), t>0$ $P-$almost surely.
We shall denote $P(\ccdot|\,x)$ the distribution of the solution
to the Cauchy problem with initial data equal to $x\in C^-_\rho$.

A theorem providing global existence and uniqueness of solutions to the Cauchy
problem for the case where the drift coefficient is
locally Lipschitz with respect to $\|\cdot\|_\rho$ is given in
Appendix~\ref{app:Cauchy problem}.
Some other existence and uniqueness results for the Cauchy problem can be found in
\cite{b:ItoNisio}, \cite{b:Protter90}.

We clearly cannot proceed without some control on the growth of
solutions in time. We obtain the control by assuming a Lyapunov--Foster
structure in the problem. As we are concerned with equations with
memory, we allow the Lyapunov functions to have memory.

\remark[Lyapunov Function]We will call a function
$V:C^-\to\R\cup \{+\infty\}$ a \textsc{Lyapunov function} for equation
\eqref{sdde} if
\begin{enumerate}
\item $V(x)\ge C_0 |x(0)|^l$   for some $C_0,l>0$ and all $x\in C^-$.
\item If a measure $P$ defines a solution $(X,W)$ of \eqref{sdde}
on
$[T_1,T_2]\subset\R$ and $P\{V(\pi_{T_1}X)<+\infty\}=1$ then
$P\{V(\pi_{t}X)<+\infty, t\in[T_1,T_2] \}=1$
  and $V(\pi_{t}X)$ satisfies the following
It\^o equation on this interval:
  \begin{equation*}
    dV(\pi_tX) = h(\pi_tX)dt + f(\pi_tX) d\widetilde W(t).
  \end{equation*}
 Here $h:C^-\to\R$ is a
  function satisfying
  \begin{equation*}
    h(x)<C_1 - C_2 V(x)^\gamma
  \end{equation*}
  for some constants $C_1,C_2,\gamma>0$ and $f:C^-\to\R$ is a function
  satisfying
  \begin{equation*}
    |f(x)| \le C_3V(x)^\delta,
  \end{equation*}
  for some $\delta \in\left[0, (1+\gamma)/2\right)$ and $C_3>0$.
  Finally $\widetilde W$ is a standard one-dimensional
  Wiener process adapted to the flow generated by $(X,dW)$.
\end{enumerate}

To show that this definition is natural, consider the Markov case where
the drift coefficient $a(x)=a(x(0))$ depends only on the present of
a trajectory $x\in C^-$. The results of \cite{b:Veretennikov97} imply
that if
\begin{equation}
\langle a(x),x(0)\rangle\le C_1-C_2|x(0)|^\alpha
\label{eq:Markov Lyapunov condition}
\end{equation}
for some
positive $C_1,C_2$ and $\alpha$ then there exists a stationary solution which
is unique. But condition \eqref{eq:Markov Lyapunov condition} means that
$V(x)=x(0)^2$ is a Lyapunov function which immediately follows from
the It\^o formula:
\begin{align*}
  d \lvert X(t)\rvert ^{2}=& 2\left[ \langle a(\pi_t(X)),
    X(t)\rangle +\frac{d}{2} \right]dt+ 2\langle X(t),dW(t)\rangle \\
  \le & 2\left[C_1 +\frac{d}{2} - C_2|X(t)|^\alpha\right]dt +
  2|X(t)|\left\langle
  \frac{X(t)}{|X(t)|},dW(t)\right\rangle.
\end{align*}

If a system possesses a Lyapunov function then it does not
fluctuate very strongly along  typical
trajectories. In
addition $V(\pi_tX)^\gamma$ averages with respect to time $t$ to a
value less than $\frac{C_1}{C_2}$
and the fluctuations of the average are not strong. To make this more
precise we define for $x \in C^-$ the fluctuations of the Lyapunov
function, denoted $\mathcal{F}V(x,t)$, by
\begin{equation*}
  \mathcal{F}V(x,t)= \left| \int_0^t V(\pi_s x)^\gamma ds \right|
  - \frac{C_1}{C_2}|t|.
\end{equation*}
For any $\rho >0$ we define the set of \textsc{nice paths} $\Nice$ by
\begin{equation*}
  \Nice=\left\{ x \in C^- : \limsup_{t\leq 0}
  \frac{|V(\pi_tx)|+ |\mathcal{F}V(x,t)| }{1+|t|^\rho} < \infty \right\}
\end{equation*}
By the first property of a Lyapunov function $\Nice \subset C^-_{\rho/l}$.
We shall say that $x\in C^-$ \textsc{averages well} if
\begin{equation*}
\limsup_{T\to\infty}\frac1T\int_{-T}^0V(\pi_t x)^\gamma dt \le
\frac{C_1}{C_2}  \ .
\end{equation*}
Notice that if $\rho <1$ then all of the path in $\Nice$ average well.
Lastly we will say that a function $g:C^-\rightarrow \Rd$ is \textsc{locally Lipschitz
on $C_\rho^-\cap \Nice[r]$} if it satisfies the Lipschitz
condition with respect to the norm $\|\ccdot\|_\rho$ on the set
\begin{equation*}
  \left\{ x \in C^- :  \limsup_{t\leq 0}\frac{|x(s)|}{1+|t|^\rho}+
  \frac{|V(\pi_tx)|+ |\mathcal{F}V(x,t)|}{1+|t|^r} < K   \right\} \ .
\end{equation*}
for any $K >0$.

For any function $F$ on $C^-$, denote $\mathcal{D}(F)=\{ x \in C^-: F(x) <
\infty\}$. We will often speak of a $\DV$-valued solution $P$ on
a set $\mathcal{R}\subset\R$. By this
we mean, $P\{ \pi_tX \in \DV, t\in\mathcal{R} \}=1$.

If $x\in C^-$, $y\in C^+$, and $x(0)=y(0)$ we define $x \cat y\in C$ as the 
concatenation of $x$ and $y$:
$$
x \cat y (t) = \left\{
\begin{array}{ll} 
x(t),& t<0,\\
y(t),& t\ge 0.
\end{array}
\right.
$$

We now state our main existence and uniqueness theorems for stationary
solutions. In Section \ref{sec:examples}, we give a few concrete
systems where the central theorems of the paper apply.  Almost all
results in this paper deal only with $\DV$-valued solutions; with the
exception of Theorem \ref{thm:uniqueness}, we say nothing about the
existence and uniqueness of other solutions.
\begin{theorem}
  \label{existence through Lyapunov}
  Let the SDE \eqref{sdde} admit a Lyapunov function $V$ and suppose
  there is an $x_0\in C^-$ such that $\bar V(x_0)=\sup_{t\in\R_-}V(\pi_tx_0)<\infty$.
  If there
  exists a finite Borel measure $\nu$ defined on subsets of $\R_-$ and
  constants $\beta, K>0$ such that
  \begin{equation}
    |a(x)| \le K +  \int_{\R_-}|V(\pi_s x)|^\beta \nu(ds),\quad x\in C^-,
    \label{measure of dependence Lyapunov case}
  \end{equation}
  the drift coefficient $a(\ccdot)$ is locally Lipschitz on
  $C_\rho^-\cap \Nice[r]$ for some $\rho > 0$,  and $r > \frac12$,
  then there exist a probability measure $P$
  on the space $\Omega$ which defines a stationary $\DV$-valued
  solution of equation {\rm \eqref{sdde}}.
\end{theorem}

\begin{theorem}
\label{thm:uniqueness} Consider a sequence of sets
$(\A_{n})_{n\in\N}$ such that $\A_n\subset \A_{n+1}\subset C^-$ and a
sequence of sets $(\B_{n})_{n\in\N}$ such that $\B_n\subset
\B_{n+1}\subset C^+$. Denote $\A_\infty=\bigcup_{n=1}^\infty A_n$,
$\B_\infty=\bigcup_{n=1}^\infty B_n$ and suppose that these sequences
and a set $\mathcal{G}\subset C$
satisfy the following properties:
\begin{enumerate}
 \item For any measure $Q$ which defines a $\mathcal{G}$-valued
 stationary solution
 to equation \eqref{sdde}
 $$Q\{\pi_0 X\in\A_\infty\}=1\mbox{\ \rm and\ }Q\{\pi_+ X\in\B_\infty\}>0.$$
 \item \label{i:contraction} For any $n\in\N$ there exists a positive function
 $\mathcal{K}_n:\R_+\to\R_+$ such that
 $\int_{\R_+}\mathcal{K}_n(t)dt<\infty$
 and for any pair of trajectories  $x_1,x_2 \in \A_n$
 and a trajectory $y\in\B_n$
 with $x_1(0)=x_2(0)=y(0)$
 \begin{equation*}
      \big\lvert a \big(\pi_t (x_1 \cat y) \big) - a \big(\pi_t
      (x_2 \cat y)
      \big)  \big\rvert ^2 \leq  \mathcal{K}_n(t).
 \end{equation*}
\end{enumerate}
Then there is at most one
$\mathcal{G}$-valued stationary solution.
\end{theorem}

\remark { The sequence of sets $A_n$ and $B_n$ can be replaced with a
  single pair of sets $A$ and $B$ where  one only requires
  $Q\{\pi_0 X\in\A, \pi_+ X\in\B\}>0$. However one must
  add the additional requirement that for any two stationary measure
  $Q_1$ and $Q_2$, $Q_1\{X(0)\in\ccdot \mbox{ and }
  \pi_0 X \in\mathcal{A},  \pi_+ X\in\B\}$ is not singular relative to
  $Q_2\{X(0)\in\ccdot \mbox{ and } \pi_0 X \in\mathcal{A}, \pi_+ X\in\B\}$.

The following theorem is a corollary of the proceeding one, but
provides simpler to verify conditions which cover many settings.
\begin{theorem}\label{thm:uniquenessPractical}
  Let the SDE \eqref{sdde} admit a Lyapunov function $V$. Fixing
  $\rho >0$ and $r > \frac12$, for $n \in \N$ define the
  sets
  \begin{align*}
    \mathcal{A}_n(\rho,r)=&\Big\{ x \in C^- : \limsup_{t\leq 0}
      \frac{|x(t)|}{1+|t|^\rho} +
      \frac{|V(\pi_tx)|+ |\mathcal{F}V(x,t)|}{1+|t|^r} < n     \Big\}\\
    \mathcal{B}_n(\rho,r)=&\Big\{ y \in C^+ :\ \mbox{\rm for all\ }x\in\mathcal{A}_n(\rho,r)
    \\ &\qquad\limsup_{t\geq 0}
      \frac{|y(t)|}{1+|t|^\rho} + \frac{|V(\pi_t(x\cat y))|+
        |\mathcal{F}V(x\cat y,t)|}{1+|t|^r} < n \Big\} \ .
  \end{align*}
  If condition \ref{i:contraction}) of Theorem \ref{thm:uniqueness}
  holds with these families of sets then there exists at most one
  $\DV$-valued stationary solution.
\end{theorem}

\section{Existence of Stationary Solutions}

The proof of the existence of stationary solutions will proceed
through a weak-limit point argument applied to the
Krylov--Bogoljubov measures (see
\cite{b:Sinai94},\cite{b:ItoNisio}).

\bpf[Proof of Theorem \ref{existence through Lyapunov}] Let $P_0$
denote a law defining a solution for the Cauchy problem for the
initial data ${x_0}$ and $P_s$ denote the time $s$-shift of this distribution
i.e. a solution of the Cauchy problem subject to initial data
$x_0(t  - s)$
defined for  $t \in (-\infty,s],s\in\R$. Formally
$P_s=P_0\theta_s^{-1}$ where $\theta_s(f,g) = (\widetilde f,\widetilde
g)$, $\widetilde f(t) = f(t-s)$, $\widetilde{ g}(t) = g(t-s)-g(-s)$.

Since the function $P_s(E)$ is  measurable with respect to $s$ for all $E\in\B$
(see \cite{b:ItoNisio}),
for $T>0$ one can define a probability measure
\begin{equation*}
  Q_T(\ccdot)=\frac{1}T\int_{-T}^0 P_s(\ccdot)ds
\end{equation*}
on the space $(\Omega,\B)$.
We will show that the family of measures $\{Q_T\}$ is tight
on a subset of $\Omega$ in an appropriate topology, which
we shall describe now.

Denote $C_\rho$ the set of trajectories $x$ in $C$ such that
$\pi_0 x\in C^-_\rho$. Define the following metric on $C_\rho$:
\begin{equation*}
  d_\rho(x,y)=\sum_{n=1}^\infty 2^{-n}(1\wedge\|\pi_n x -\pi_n y\|_\rho).
\end{equation*}
Finally equip the space
$\Omega_\rho=C_\rho\times C$ with product topology of
$d_\rho-$topology in $C_\rho$ and LU-topology in $C$.

The general idea of proof of tightness in a space of continuous trajectories
is to obtain uniform bounds for marginal distributions and for distributions
of increments of the trajectories, see \cite{b:Billingsley68}.

\begin{lemma}
  \label{l:QTmoments}
  For any $\kappa\ge0$, $T\geq 1$ and $S\ge0$ the moments $\E_{Q_T} V(\pi_t X)^\kappa$
  are uniformly bounded for $t\in(-\infty,S]$.
\end{lemma}
\bpf[Proof of Lemma \ref{l:QTmoments}]
In the sequel, we will write $V_t$ instead of $V(\pi_t X)$ for
brevity.  We may here assume without loss of generality that
$V(\ccdot)\ge1$ (indeed, if $V$ is a Lyapunov function so is $V+1$ with
possibly different choice of constants $C_1$ and $C_2$)
and hence the moments of negative order are uniformly bounded.

We prove first that for each positive $\kappa$
there exists a finite constant $N_\kappa$ such that
for any $T>1$ and any $t\in[-T,0]$
\begin{equation}
\label{eq: uniform in t,T}
  \frac1T\E_{P_0}\int_0^{T+t}V_s^\kappa\le N_\kappa.
\end{equation}

Define $\tau_R(t) = t\wedge\tau_R$ where $\tau_R=\inf\{t:V_t\ge R\}$.
Let $m\ge1$. The It\^o formula and the assumptions on Lyapunov function
$V$ imply that for $T>0$
\begin{multline*}
  V_{\tau_R(T+t)}^{m}-V_0^{m}\le m\int_0^{\tau_R(T+t)}
  \Biggl\{V_s^{m-1}[C_1-C_2V_s^\gamma] \\ +
    \frac{(m-1)}{2}C_3^2V_s^{m-2+2\delta}\Biggr\}ds+
  m\int_0^{\tau_R(T+t)} V_s^{m-1}f(\pi_tX) d\widetilde W(s).
\end{multline*}
holds $P_0$-a.s. Taking expectations of both sides, passing to limit
$R\to\infty$, using the regularity of the solution ($\tau_R\to\infty$
as $R\to\infty$) established in the Appendix \ref{app:Cauchy problem}
we get
\begin{multline*}
  \E_{P_0}\left[V_{T+t}^{m}  -  V_0^{m}\right]
  \le mC_1\E_{P_0}\int_0^{T+t} V_s^{m-1}ds \\ -
  mC_2\E_{P_0}\int_0^{T+t} V_s^{m-1+\gamma}ds
  +\frac{m(m-1)}{2}C_3^2\E_{P_0}\int_0^{T+t}V_s^{m-2+2\delta}ds.
\end{multline*}

Dividing both sides of this inequality by $T$ and using positivity
of $V_{T+t}$ we obtain that
\begin{multline*}
 \frac{mC_2}{T}\E_{P_0}\int_0^{T+t} V_s^{m-1+\gamma}ds\le
 \frac{mC_1}{T}\E_{P_0}\int_0^{T+t} V_s^{m-1}ds
 \\ +\frac{m(m-1)}{2T}C_3^2\E_{P_0}\int_0^{T+t}V_s^{m-2+2\delta}ds
 +\frac{V({x_0})^m}{T}.
\end{multline*}

Since
$\eps=[1+\gamma-2\delta]\wedge\gamma>0$ we can apply the last inequality
iteratively thus  extending the domain of applicability of
\eqref{eq: uniform in t,T}  at each iteration by $\eps$.

There are three cases to consider, namely $t\in(-\infty,-T)$, $t \in
[-T,0)$ and $t \in [0,S)$. Suppose now that $t\in[-T,0]$. Then
\begin{equation*}
  \E_{Q_T}V_t^\kappa=\frac{1}{T}\E_{P_0}\int_t^{t+T}V_s^\kappa ds
  \\=\frac1T\E_{P_0}\int_0^{t+T}V_s^\kappa ds+
  \frac{t}{T}\bar V(x_0)^\kappa\le N_\kappa+\bar V({x_0})^\kappa
\end{equation*}
and the same estimate is obviously true if $t<-T$.

Suppose now $t\in[0,S]$. Then
\begin{multline*}
\E_{Q_T} V_t^\kappa=
\frac{1}{T}\E_{P_0}\int_{t}^{t+T} V_s^\kappa ds\le
\frac{T+S}{T}\E_{Q_T}V_0^\kappa
\le
\frac{T+S}{T}[N_\kappa+\bar V(x_0)^\kappa].
\end{multline*}
So for any $\kappa$ we have proved that $\E_{Q_T} V_t^\kappa$
is bounded uniformly with respect to $T>1$ and $t\in(-\infty,S]$.
\epf

Property 1 of Lyapunov function $V$ and Lemma \ref{l:QTmoments}
imply that all the moments of $X$ are also uniformly bounded.

Now let us estimate increments of the process $X$.
\begin{lemma}\label{l:increments}
  There exists a constant $C>0$ so that for any $t_1$,
  $t_2$ and $T\geq 1$, one has $ Q_T\{|X(t_2)-X(t_1)|> z\}\le C(z^{-4} + z^{-2})|t_2-t_1|^2$.
\end{lemma}
\bpf[Proof of Lemma \ref{l:increments}]
\begin{align}
  \label{increments Lyapunov case}
  Q_T\{|X(t_2)-&X(t_1)|> z\} \notag\\& \leq Q_T\{|W(t_2)-W(t_1)|>z/2\}+
  Q_T\Big\{\int_{t_1}^{t_2} a(\pi_{\theta}X)d\theta >z/2 \Big\} \notag \\
  & \leq\frac{16}{z^{4}}\E_{Q_T} |W(t_2)-W(t_1)|^4 +
  \frac{4}{z^{2}}\E_{Q_T}\Big(\int_{t_1}^{t_2}
  a(\pi_{\theta}X)d\theta\Big)^2 .
\end{align}

The first term can be estimated through the well-known expression for
moments of Gaussian distribution.  To estimate the second term we use
the Fubini theorem, elementary inequality $|xy|\le(x^2+y^2)/2$ and
relation \eqref{measure of dependence Lyapunov case}:
\begin{equation}
  \label{increments due to drift Lyapunov case}
  \begin{split}
    \E\int_{t_1}^{t_2}&a\big(\pi_\theta(x)\big)d\theta =
    \int_{t_1}^{t_2}\int_{t_1}^{t_2}\E
    a(\pi_{\theta_1}(x))a\big(\pi_{\theta_2}(x)\big) d\theta_1d\theta_2
    \\ &\le (t_2-t_1)^2 \sup \E a\big(\pi_\theta(x)\big)^2\\&\le
    (t_2-t_1)^2 \sup\E\Big(K+\int_{\R_-}V(t+s)^\beta \nu(ds)\Big)^2\\
    &\le (t_2-t_1)^2 \Big(K^2+ 2K\int_{\R_-}\sup\E V(t+s)^\beta
    \nu(ds)\\& \qquad \qquad \qquad + \int_{\R_-}\int_{\R_-}\sup\E
    V(t+s_1)^\beta V(t+s_2)^\beta \nu(ds_1) \nu(ds_2) \Big) \\ &\le
    (t_2-t_1)^2 \left(K^2 + 2K(M+1)^\beta\nu(\R_-)+ (M+1)^{2\beta}
      \nu(\R_-)^2 \right).
  \end{split}
\end{equation}
Inequalities \eqref{increments Lyapunov case} and \eqref{increments
  due to drift Lyapunov case} imply
\begin{equation*}
  Q_T\{|X(t_2)-X(t_1)|> z\}\le 48z^{-4}|t_2-t_1|^2 + Cz^{-2}|t_2-t_1|^2.
\end{equation*}
for a constant $C>0$.\epf

Since we have uniform moment estimates and uniform increment estimates
in probability, the tightness of the distribution of the process $X$,
in the uniform topology on any finite interval under the measures
$Q_T$, follows immediately from [3, Theorem 12.3] \cite[Theorem
12.3]{b:Billingsley68}.

To finish the proof of tightness in $\Omega_\rho$ we need the following lemma.
\begin{lemma}\label{lem:Krylov--Bogolyubov envelope}
For any $\sigma>0$ and any $t\in\R$ the random variable
$\|\pi_t X\|_\sigma$ is finite $Q_T$-a.s. for any $T>0$. Moreover,
the family of distributions $Q_T\{\|\pi_t X\|_\sigma\in\ccdot\}$ is
tight.
\end{lemma}
\bpf[Proof of Lemma \ref{lem:Krylov--Bogolyubov envelope}] We omit the
proof which is similar to the proof of Lemma \ref{l:fixPointNice}
below and Theorem \ref{asymtotic behavior} from the next section. It
relies on uniform estimates of marginals and increments in probability
and the Borel--Cantelli lemma.  \epf

Since the distribution of $W$ in $C$ is the same under all measures
$Q_T$, it suffices to demonstrate tightness of the distributions
of $X$. Consider any $S>0$ and fix $\eps>0$. Choose any $\sigma<\rho$
and use Lemma \ref{lem:Krylov--Bogolyubov envelope} to find  $K_\eps>0$
such that
\begin{equation*}
  Q_T\big\{\|\pi_t X\|_\sigma>K_\eps\big\}>1-\frac\eps2.
\end{equation*}
For any $n\in N$  due to tightness of distribution of $X$ on $[-n,S]$ one
can choose a compact set $E_n\subset C_{[-n,S]}$ such that
\begin{equation*}
  Q_T\{X[-n,S]\in E_n\}>1-2^{-n-1}\eps.
\end{equation*}
Let $E_\infty=\big\{x\in C_{(-\infty,S]}\;\big|\; x[-n,S]\in E_n,
n\in\N\ \mbox{and\ }\|x\|_\sigma\le K_\eps \big\}$.

Since $\sigma<\rho$ it is straightforward to show
that $E_\infty$ is compact in the norm
$\|\pi_S\ccdot\|_\rho$ and $Q_T\{\pi_S X\in E_\infty\}>1-\eps$.

So for any $S\in\N$ and $\eps$ we can build a set $E^{(S)}$ which is
compact in the corresponding norm and $Q_T\{X(-\infty,S]\in E^{(S)}\}>1-\eps2^{-S}$.
Use the same construction to build a compact set $E\subset C_\rho$ with
$Q_T\{X\in E\}>1-\eps$.

So the tightness of the distributions of $X$ in $C_\rho$ under $Q_T$
is proved.
The classical Prokhorov theorem implies that
$Q_{T_n}\stackrel{Law}{\to} Q_\infty$ when $n\to\infty$ for some
infinitely increasing sequence $(T_n)_{n\in\N}$.

To conclude the proof of Theorem \ref{existence through Lyapunov} we
need to establish some properties of the trajectories on which
$Q_\infty$ is concentrated.
\begin{lemma}
  \label{l:fixPointNice} For any $\varkappa>\frac12$ and $\rho>0$
 \begin{equation*}
Q_\infty\left\{\sup_{t \leq 0 } \frac{|\mathcal{F}V(X,t)|}{1+|t|^\varkappa} <
  \infty\right\}=Q_\infty\left\{\sup_{t \leq 0 } \frac{V(\pi_tX)}{1+|t|^\rho} <
  \infty\right\}=1 \ .
\end{equation*}
\end{lemma}
\bpf[Proof of Lemma \ref{l:fixPointNice}]

First notice that due to the construction of $Q_T$ as an average over
initial value problems, if $T$ be chosen so that $n+1<T2^{-n}$
then for any set $A\in\Omega$ the probability $Q_{T}(A)$ is bounded
by $Q_{T}(A\cap E)+2^{-n}$ where for any $\omega\in E$ the equation
\eqref{sdde} is satisfied on $[-n,0]$. This will allow us to estimate
$Q_T$ using the dynamics up to a small error.

We begin with the first claim controlling $\mathcal{F}V(\pi_tX)$. The
claim is implied   by
\begin{equation*}
  \sum_{n=1}^\infty Q_\infty\left\{\int_{-S}^0V_s^\gamma ds
    - \frac{C_1}{C_2}S>3n^\varkappa, S\in[-n-1,-n]\right\}<\infty
\end{equation*}
which in turn follows from
\begin{equation}
  \limsup_{m\to\infty} Q_{T_m}\left\{\int_{-S}^0V_s^\gamma ds
    - \frac{C_1}{C_2}S>3n^\varkappa, S\in[-n-1,-n]\right\}\leq q_n
  \label{eq:majorization}
\end{equation}
and $\sum_nq_n<\infty$.

Let $T_m$ be chosen so that $n+1<T_m2^{-n}$ then as mentioned above for any set
$A\in\Omega$ the probability $Q_{T_m}(A)$ is bounded by $Q_{T_m}(A\cap
E)+2^{-n}$ where for any $\omega\in E$ the equation \eqref{sdde} is
satisfied on $[-n,0]$ and hence on $E$ one has
\begin{equation*}
  V_0-V_{-S}\le\int_{-S}^0 (C_1 - C_2V_t^\gamma)dt+\int_{-S}^0f(\pi_t X)d\widetilde{W}(t).
\end{equation*}
So,
\begin{equation*}
  \int_{-n-1}^0  V_t^\gamma dt- \frac{C_1}{C_2}{(n+1)}\le
  \frac{V_{-n-1}-V_0}{C_2}+\frac1{C_2}\int_{-n-1}^0f(\pi_t X)d\widetilde{W}(t)
\end{equation*}
and
\begin{multline}
\label{eq:averaging-martingale part}
  Q_{T_m}\left\{\int_{-n-1}^0  V_t^\gamma dt- \frac{C_1}{C_2}>2n^\varkappa\right\}
  \le Q_{T_m}\{V_{-n-1}-V_0>C_2n^\varkappa\}\\
   +  Q_{T_m}\left\{\int_{-n-1}^0f(\pi_t X)d\widetilde{W}(t)> C_2n^\varkappa\right\}
   \le C(p)n^{p\left(\frac12-\varkappa\right)}\le C(p)n^{-2}
\end{multline}
for $\varkappa>1/2$ and  large enough $p$
where the last line follows from uniform boundedness of all moments of $V$ and
Burkholder's inequality (see \cite[Theorem 54]{b:Protter90}).
\begin{multline}
\label{eq:averaging-time1increment}
  Q_{T_m}\left\{\sup_{S\in[-n-1,-n]}\left[\int_{-n-1}^SV_s^\gamma ds
  -\frac{C_1}{C_2}(S+n+1)\right]>n^\varkappa\right\}
  \\ \le Q_{T_m}\left\{\int_{-n-1}^{-n}V_s^\gamma ds >n^\varkappa\right\}\le C n^{-2}.
\end{multline}
where we used boundedness of all moments of $V$ in the last estimate.

Now \eqref{eq:majorization} with $q_n=Cn^{-2}+2^{-n}$ follows from
\eqref{eq:averaging-martingale part} and
\eqref{eq:averaging-time1increment} and control of the time average
fluctuations claimed by Lemma \ref{l:fixPointNice} is proved.

We now turn to the remaining claim controlling $V(\pi_t X)$. Choosing
$T_m$ as above,  for $-n\in\N$
\begin{multline}
\label{supA}
  Q_{T_m}\left\{\sup_{t\in[n,n+1]} V_t> K(|n|^\rho+1)\right\}
  \le Q_{T_m}\left\{V(n)>\frac{K}{2}(|n|^\rho+1)\right\}\\+
  Q_{T_m}\left\{\sup_{t\in[n,n+1]}V_t-V(n)>\frac{K}{2}(|n|^\rho+1)\right\}=I_1+I_2.
\end{multline}
Next we estimate both terms using Chebyshev's inequality and the
uniform bounds from Lemma \ref{l:QTmoments}.
\begin{equation*}
  I_1\le\frac{\E_{Q_{T_m}} V(n)^\kappa}{\frac{K}{2}(|n|^\rho+1)^\kappa}\le\frac{C}{n^2+1}
\end{equation*}
for some constant $C>0$ if $\kappa>2/\rho$.
\begin{multline}
  I_2\le Q_{T_m}\left\{\sup_{t\in[n,n+1]}\int_n^tf(\pi_sX)d\widetilde W(s)
    > \frac{K}{2}(|n|^\rho+1)-C_1\right\}\\
  \le \frac{C\E_{Q_{T_m}}
    \left[\sup_{t\in[n,n+1]}\int_n^tf(\pi_sX)d\widetilde
      W(s)\right]^{2p}}{n^2+1}
  \label{I_2A}
\end{multline}
for some constant $C>0$ and  $p\in\N$ such that $2p\rho>2$.

To prove that the expectation in the right-hand side of \eqref{I_2A} is
finite use Burkholder's inequality and the uniform estimates on
$E_{Q_{T_m}} V_t^p$ given by Lemma \ref{l:QTmoments}:
\begin{align}
  \label{eq:BDGBound}
  \E_{Q_{T_m}}&\left[\sup_{t\in[n,n+1]}\int_n^tf(\pi_sX)d\widetilde
    W(s)\right]^{2p}\le
  K_{2p}\E_{Q_{T_m}}\left[\int_n^{n+1}f^2(\pi_sX)ds\right]^{p} \\
    \notag &=
  K_{2p}\E_{Q_{T_m}}\int_n^{n+1}\ldots\int_n^{n+1}f^2(\pi_{s_1}X)\ldots
  f^2(\pi_{s_1}X)ds_1\ldots ds_p \notag \\ &\le
  K_{2p}\E_{Q_{T_m}}\int_n^{n+1}\ldots\int_n^{n+1}\left[f^{2p}(\pi_{s_1}X)+\ldots
    + f^{2p}(\pi_{s_1}X)\right]ds_1\ldots ds_p \notag \\ &\le K_{2p}p
  C <\infty \ . \notag
\end{align}
To complete the proof apply
\eqref{supA}--\eqref{I_2A} and the Borel--Cantelli lemma as in the first
part.  \epf

\bpf[Completion of Proof of Theorem \ref{existence through Lyapunov}]
All that remains in the proof of Theorem \ref{existence through
  Lyapunov} is to show that $Q_\infty$ is concentrated on solutions
which solve the equation. From Lemma \ref{l:fixPointNice}, we see that
for any $\rho >0$ and $r>\frac12$, $Q_\infty\{ X \in C^\rho_-\cap
\Nice[r]\} =1$.  Since the drift coefficient $a(\ccdot)$ is
continuous on the set of such paths, the reasoning from
\cite[p.21--25]{b:ItoNisio} shows that $Q_\infty$ defines a stationary
solution of the equation \eqref{sdde}. Theorem \ref{existence through
  Lyapunov} is proved.  \epf

\section{Properties of Stationary Solutions.}
Before turning to the question of what additional requirements are
sufficient to guarantee the uniqueness of the stationary measure, we
extract a number of important properties which any stationary measure
must possess given the assumptions already made. Specifically, we give
bounds on the moments and growth to the Lyapunov function in time,
prove that the averaging property is fulfilled a.s.  and characterize
the marginals of any stationary measure at a fixed given time.

\subsection{Control of Moments and Asymptotic Path Behavior}
\label{section moments}

\begin{theorem}
  \label{momestimates} Under the conditions of Theorem
  \ref{existence through Lyapunov}, for any $\kappa \in \RR$ there
  exist a single, fixed constant $M_k$ so that
  \begin{equation*}
    \E_Q V_t^\kappa \leq M_\kappa <\infty
  \end{equation*}
  under any measure $Q$ which defines a stationary $\DV$-valued
  solution $X$ of equation \eqref{sdde}.
\end{theorem}

\bpf Let $g_{N,m}(x)=(x\wedge N)^m$ for $x\in\R, m, N>0$.  Apply
the It\^o-Meyer formula (see \cite[Theorem 51]{b:Protter90}) to
$g_{N,m}(V_t)$:
\begin{multline*}
  g_{N,m}(V((T)) - g_{N,m}(V((T)) \\ \le\int_0^{T}\ONE\{V_t\le N\}
 m\Bigl[V_t^{m-1}(C_1-C_2V_t^\gamma) +
  \frac{m-1}{2}C_3^2V_t^{m-2+2\delta}\Bigr]dt+
\\
  m\int_0^{T} V_t^{m-1}f(\pi_tX)\ONE\{V_t\le N\} d\widetilde W(t)-\psi(T).
\end{multline*}
Here $\psi$ is a non-decreasing function such that $\psi(0)=0.$

Take expectations of both sides of the last inequality with respect
to the stationary measure $Q$:
\begin{multline*}
  C_2m\E_Q\int_0^{T} V_t^{m-1+\gamma}\ONE\{V_t\le N\}dt\\ \le
  C_1m\E_Q\int_0^{T}V_t^{m-1}\ONE\{V_t\le N\}dt \\ +
  \frac{m(m-1)}{2}C_3^2\E_Q\int_0^{T}V_t^{m-2+2\delta}\ONE\{V_t\le
  N\}dt.
\end{multline*}
Take the limit $N\to\infty$ and use stationarity of $Q$ to get
\begin{equation*}
  \E_{Q}V_t^{m-1+\gamma} \le \frac{C_1}{C_2}
  \E_{Q}V_t^{m-1} +
  \frac{(m-1)C_3^2}{2C_2} \E_{Q}V_t^{m-2+2\delta}.
\end{equation*}

As in the previous section since $\eps=1+\gamma-2\delta>0$ we can
apply the last moment inequality iteratively to see that moments of
$V_t$ are bounded under the stationary measure $Q$. The proof is
complete.  \epf

\begin{theorem}
\label{asymtotic behavior}  Let the SDE \eqref{sdde} admit a Lyapunov
function $V$. Suppose measure $Q$ defines a stationary $\DV$-valued
solution for \eqref{sdde}. Then for any $\rho>0$
\begin{equation*}
  Q\left\{\sup_t\frac{V(\pi_t X)}{1+|t|^\rho}< \infty \right\}=1.
\end{equation*}
\end{theorem}
\bpf We proceed as in the proof of the second claim in Lemma
\ref{l:fixPointNice}. As in \eqref{supA}, for $n\in\Z$
\begin{multline}
\label{sup}
Q\left\{\sup_{t\in[n,n+1]} V_t> K(|n|^\rho+1)\right\} \le
Q\left\{V(n)>\frac{K}{2}(|n|^\rho+1)\right\}\\+
Q\left\{\sup_{t\in[n,n+1]}V_t-V(n)>\frac{K}{2}(|n|^\rho+1)\right\}=I_1+I_2.
\end{multline}
Next we estimate both terms using Chebyshev's inequality. Namely,
\begin{equation}
  I_1\le\frac{\E_Q V(n)^\kappa}{\frac{K}{2}(|n|^\rho+1)^\kappa}\le\frac{C}{n^2+1}
\label{I_1}
\end{equation}
for some constant $C>0$ if $\kappa>2/\rho$. Calculations analogous to
\eqref{I_2A} and \eqref{eq:BDGBound} give
\begin{equation}
  I_2\le \frac{ K_{2p}p
  \E_Qf^{2p}(\pi_TX)}{n^2+1}\le \frac{K_{2p}p C_3^{2p}\E_Q V(T)^{2p\delta}}{n^2+1}
  \label{I_2}
\end{equation}
for some constant $C>0$ and $p\in\N$, such that $2p\rho>2$, and for
any time $T \in \RR$. (By stationarity the choice of $T$ does not
matter.) To complete the proof apply \eqref{sup}--\eqref{I_2} and the
Borel--Cantelli lemma.  \epf.

\begin{theorem}\label{thm:avgFluc}
  Let the SDE \eqref{sdde} admit a Lyapunov function $V$. Suppose $Q$
  is a measure which defines a stationary $\DV$-valued solution to
  equation \eqref{sdde}. Then for any $\varkappa > \frac12$
\begin{equation*}
Q\left\{\sup_{t}\frac{|\mathcal{F}V(X,t)|}{1+|t|^\varkappa} <
  \infty\right\}=1  \ .
\end{equation*}
\end{theorem}
\bpf The needed calculations parallel those in the proof of the second
part of Lemma \ref{l:fixPointNice}; we give most of the details
nonetheless. Using the definition of Lyapunov function we have
\begin{align*}
  V_0+ \int_{-T}^0 C_2V_t^\gamma dt- C_1|T| &\leq V_{-T} +
  \int_{-T}^0f(\pi_t X)d\widetilde{W}(t)\\
  V_T+ \int^{T}_0 C_2V_t^\gamma dt- C_1|T| &\leq V_{0} +
  \int^{T}_0 f(\pi_t X)d\widetilde{W}(t)
\end{align*}
Since by Theorem \ref{asymtotic behavior}, $V_{-T}$ grows slower than
$T^\frac12$ with probability 1 it suffices to show that the last term
in each of the above lines grows no faster than $T^\varkappa$.
Denote $I(t_1,t_2)=\int_{t_1}^{t_2}f(\pi_t X)d\widetilde{W}(t)$ for
$t_1<t_2$.  We shall prove that
\begin{multline}
  \sum_{n=1}^\infty Q\biggl\{|I(-n,0)|>n^{\varkappa} \biggr\} +
  Q\biggl\{\sup_{s\in[-n,-n+1]}|I(-n,s)|>n^{\varkappa} \biggr\}\\ +
  \sum_{n=1}^\infty Q\biggl\{ \sup_{s\in[0,n+1]}|I(0,s)| >
    2n^\varkappa \biggr\}<\infty.
\label{eq:growth of martingale part}
\end{multline}
Since the sum
\begin{equation*}
  \sum_{n=1}^\infty Q\left\{\sup_{s\in[-n,-n+1]}
  |I(s,0)|>2n^{\varkappa} \right\}+ \sum_{n=0}^\infty Q\left\{\sup_{s\in[n,n+1]}
  |I(0,s)|>2n^{\varkappa} \right\}
\end{equation*}
is majorized by by the previous one, the applying the Borel--Cantelli
lemma will conclude the proof. To derive \eqref{eq:growth of
  martingale part} use Burkholder's inequality and the fact that
moments of $V$ are uniformly bounded:
\begin{equation*}
Q\left\{|I(-n,0)|>n^{\varkappa}
  \right\}\le n^{-2p\varkappa} C \E_Q\left[\int_{-n}^0V_t^{2\delta}dt\right]^p\le
  Cn^{-2p\varkappa}n^p\le C n^{p(1-2\varkappa)}
\end{equation*}
and
\begin{align*}
Q\left\{\sup_{s\in[-n,-n+1]}|I(-n,s)|>n^{\varkappa}
  \right\}&\le n^{-2p\varkappa} C \E_Q\left[\int_{-n}^{-n+1}V_t^{2\delta}dt\right]^p\le
  Cn^{-2p\varkappa} \\
 Q\left\{\sup_{s\in[n,n+1]}|I(0,s)|>n^{\varkappa}
  \right\}&\le n^{-2p\varkappa} C \E_Q\left[\int_{0}^{n}V_t^{2\delta}dt\right]^p\le
  Cn^{p(1-2\varkappa)}
\end{align*}
where by $C$ we denote possibly different constants.  Since $\varkappa
>1/2$, \eqref{eq:growth of martingale part} is finite as claimed and
the theorem is proved.  \epf

\subsection{Regularity of Time $t$ Marginals}
For any $x\in C^-$, let $P_t(\ccdot|\,x)$ denote the measure induced on
$\RR^d$ at time $t$ by the dynamics starting from the past $x$ at time
$0$. Similarly for any stationary measure $Q$ denote the marginal at
time $t$ on $\RR^d$ by $Q_t(A)=Q\{ X(t) \in A \}$. By
stationarity the measure $Q_t$ is independent of $t$. The following
Theorem along with Lemma \ref{l:infFutures} in the next section are the
key elements in the uniqueness proof.
\begin{theorem}\label{l:timeTmarginal}  Let $Q$ be any
measure defining a stationary solution of equation \eqref{sdde}.
For $Q$-almost every $x \in C^-$ and every $t>0$, $P_t(\ccdot|\,x)$ is
  equivalent to the Lebesgue measure on $\RR^d$. In addition, $Q_s$ is
  equivalent to the Lebesgue measure  on $\RR^d$ for any $s \in \RR$.
\end{theorem}
\bpf Since, by stationarity,
 $Q_t(A)=\int P_t(A\,|\,\pi_0x)Q(dx)$ the equivalence of $Q_t$ to the
 Lebesgue measure follows from
 that of $P_t$. To prove that $P_t(\ccdot\,|\,x)$ is equivalent to the
 Lebesgue measure,
 it is sufficient to show that the distribution
 $P^X_{[0,t]}(\ccdot|\,x)$ of the process $X$ on the
 interval $[0,t]$ is equivalent to the distribution
 $P^W_{[0,t]}(\ccdot|\,x)$ of standard Wiener process $W$
 started at $x(0)$.
 To apply Lemma \ref{l:compareMeasures} from the appendix,
 we need a truncation to guarantee the Novikov
 like condition in \eqref{eq:Novikov}.  We define the adapted function
 \begin{equation*}
   T_R(X)= \inf\left\{ s>0 : \int_0^s |a(\pi_s X)|^2 ds > R\right\} \ .
 \end{equation*}
  By stationarity, we know that for $Q$-almost
 every initial condition $x$,
 \begin{equation*}
   P( \exists R \mbox{ so } T_R(X) > t|\,x) = 1.
 \end{equation*}
 On the other hand for $Q$-almost every initial condition $x \in C^-$
 and any
 \begin{equation}
   P^X_{[0,t]}( \ccdot ; T_R(X) > t|\,x) \sim P^W_{[0,t]}( \ccdot ; T_R(W) > t|\,x)
 \label{eq:equivalence to wiener}
 \end{equation}
 where $\sim$ denotes equivalence of measures. Indeed, the definition
 of  the stopping
 time $T_R(X)$ guarantees the Novikov condition \eqref{eq:Novikov} and
 the equivalence in \eqref{eq:equivalence to wiener}
 is implied by
 Lemma \ref{l:compareMeasures} with $\mathcal{B}=\{T_R(X)>t\}$.
 If $R \rightarrow \infty$ then
 the sequences of measures $P^X_{[0,t]}(\ccdot ; T_{R_n}(X) > t|\,x)$
 and $ P^W_{[0,t]}(\ccdot ; T_{R_n}(W) > t|\,x)$ increase to
 $P^X_{[0,t]}(\ccdot|\,x )$ and $P^W_{[0,t]}(\ccdot |\,x )$ respectively and
 equivalence is preserved under this limit.
\epf

\section{Uniqueness of Stationary Solutions}

The proof of Theorem \ref{thm:uniqueness} rests on the following lemma
which in turn relies on the two lemmas which follow it and  contain
the heart of the matter. Given any stationary measure $Q$, for $B
\subset C^+$, we define $Q_+(B; \mathcal{B})= Q\{ \pi_+X \in B \cap
\mathcal{B})$.

\begin{lemma}
  \label{l:equivFutures} In the setting of Theorem
  \ref{thm:uniqueness}, $Q_{1+}( \ccdot ;\mathcal{B}_\infty ) \sim
  Q_{2+}( \ccdot ;\mathcal{B}_\infty )$ for any two measures $Q_{1}$ and
  $Q_{2}$ defining stationary $\mathcal{D}$-valued solutions to \eqref{sdde}.
\end{lemma}
Using this lemma, one quickly obtains a proof of  Theorem
\ref{thm:uniqueness}. Once two stationary
measures are shown not to be singular on the future then they must be
the same measure. Lemma \ref{l:equivFutures} gives the necessary
equivalence, one way to see that this ensures uniqueness is given
next.

\bpf[Proof of Theorem \ref{thm:uniqueness}]
Consider two ergodic measures  $Q_{1}$ and $Q_{2}$ defining
stationary $\Dc$-valued solutions.
Let $\phi:C \rightarrow \RR$ be an arbitrary
bounded functional which depends on values of its argument within
a finite interval.
By the Birkhoff--Khinchin ergodic theorem (see
for instance \cite{b:Sinai94}), there
exists deterministic constants $\phi_{i}$ so that
\begin{equation*}
  \lim_{T\rightarrow \infty} \frac1T \int_0^T \phi(\pi_tX)dt = \int
  \phi(x)dQ_{i}(x)=  \phi_{i}
\end{equation*}
$Q_{i}$-almost surely. Let $C^{(i)}$ be subsets of $C^+$ of full
$Q_{i}$-measure so that the above limit on the left hand side
converges to the given constant $\phi_{i}$. Since $Q_{1}\{\pi_+ X
\in C^{(1)}\} =1$, we know that $Q_{1}\{\pi_+ X \in C^{(1)} \cap
\mathcal{B}_\infty\}>0$. Then Lemma \ref{l:equivFutures} implies that
$Q_{2}\{ \pi_+X \in C^{(1)}\cap \mathcal{B}_\infty \} >0$. Since $Q_{2}\{
\pi_+X \in C^{(2)}\}=1$, we have that $Q_{2}\{ \pi_+X \in
C^{(1)}\cap C^{(2)} \cap \mathcal{B}_\infty\} >0$. Hence $C^{(1)}\cap
C^{(2)}$ is non-empty implying that $\phi_{1}= \phi_{2}$.  Since
the $\phi$ was  arbitrary, we conclude that the distribution of $X$ under
$Q_{1}$ and $Q_{2}$ is the same.
Due to the general fact that
any stationary measure can be represented as a convex combination of
ergodic ones the distribution of $X$ is the same for any
measure $Q$ defining a stationary $\Dc$-valued solution
(the ergodic decomposition of a measure defining
a $\Dc$-valued solution gives zero weight to ergodic measures defining
solutions which are not $\Dc$-valued).

Since the trajectory of $X$ on $(0,\infty)$ is fully
determined by $\pi_0X$ and the trajectory of $W$, the distribution
of $\pi_0X$ determines uniquely the joint distribution of $X$ and $W$
and the theorem is proved.
\epf

Lemma \ref{l:equivFutures} which was pivotal in the preceding proof is
itself a consequence of the following lemma which along with Theorem
\ref{l:timeTmarginal} contain the central steps in the uniqueness
proof.  We will give its proof directly after its statement, returning
to the proof of Lemma \ref{l:equivFutures} at the end of the section.

\begin{lemma}
  \label{l:infFutures}
  Suppose that assumptions of Theorem \ref{thm:uniqueness} are satisfied.
  For fixed $n\in\N$, if $x_1,x_2\in \A_n$ and
  $x_{1}(0)=x_{2}(0)$ then $P(\ccdot; \B_n|x_{1}) \sim
  P(\ccdot ; \B_n|\,x_{2})$ where
  $P(B ;\mathcal{B}_n|\,x) = P\{X\in B\cap\mathcal{B}_n|\,x\}$
  for $x\in \A_n$  and
  $B \subset C^+$.
\end{lemma}

\bpf[Proof of Lemma \ref{l:infFutures} ] We are going to derive this lemma
from Lemma \ref{l:compareMeasures} in the appendix.
Set $f_i(t,y[0,t])=a(\pi_t(x_i \cat y)), i=1,2.$
Notice that for a fixed
$t$, $f_i(t,\ccdot)$ can be thought of as functions
from $C([0,t];\RR^d)$ to $\RR^d$ since the part of each $\pi_t
(x_i \cat y)$ before 0 is fixed.

Since
\begin{multline}
  \exp \left\{\frac12 \int_0^\infty \big\lvert a \big(\pi_t( x_{1}\cat
  y) \big)
    - a \big(\pi_t(x_{2}\cat y) \big) \big\rvert ^2 dt \right\}
  \ONE_{\mathcal{B}_n}(x)
  \\ \le
  \exp \left\{\frac12\int_0^\infty
    \mathcal{K}_n(t)dt \right\}<\infty
\label{eq:bounded Novikov}
\end{multline}
condition \eqref{eq:Novikov}
is satisfied
and  Lemma \ref{l:compareMeasures} now implies the desired equivalence
of measures.
\epf

With Lemma \ref{l:infFutures} proved, we return to the proof of Lemma
\ref{l:equivFutures}.

\bpf[Proof of Lemma \ref{l:equivFutures}]
It is sufficient to show that for all $n\in\N$ if
 $B$ is a set in $C^+$ such that
$Q_{1}\{\pi_+X\in B\cap\B_n\}=0$ then $Q_{2}\{\pi_+X\in
B\cap\B_n\}=0$.

From the relation
\begin{align*}
  0=&Q_{1}\{\pi_+X\in B\cap\B_n,\pi_0X\in\A_\infty\}\\ =&
  \int_{\Rd} Q_{1}\{\pi_+X\in B\cap\B_n,\pi_0X\in\A_\infty
  |X(0)=x_0\}Q_1\{X(0)\in dx_0\}
\end{align*}
and Theorem \ref{l:timeTmarginal}
we see that for $x_0$ in some set $E\subset\Rd$ with complement
in $\Rd$ of zero Lebesgue measure
\begin{equation*}
  Q_{1}\{\pi_+X\in B\cap\B_n,\pi_0X\in\A_\infty|X(0)=x_0\}=0.
\end{equation*}
This means that for all $m\ge n$ and for
$Q_{1}\{\pi_0X\in\ccdot\cap\A_m|X(0)=x_0\}$- almost all $x\in C^-$
\begin{equation*}
  Q_{1}\{\pi_+X\in B\cap\B_n|\pi_0X=x\}=0.
\end{equation*}
Lemma \ref{l:infFutures} now implies that
if $x_0\in E$ then
for {\it all} $x\in \A_m$ with $x(0)=x_0$
\begin{equation*}
  Q_{2}\{\pi_+X\in B\cap\B_n|\pi_0X=x\}=0,
\end{equation*}
which in turn implies
\begin{equation*}
  Q_{2}\{\pi_+X\in B\cap\B_n, \pi_0X\in \A_m|X(0)=x_0\}
  =0\quad\mbox{for\ } x_0\in E
\end{equation*}
and
\begin{equation*}
  Q_{2}\{\pi_+X\in B\cap\B_n, \pi_0X\in \A_\infty|X(0)=x_0\}
  =0\quad\mbox{for\ } x_0\in E
\end{equation*}
To complete the proof
integrate the last relation over $\Rd$ with respect to
$Q_2\{X(0)\in dx_0\}$ and use Theorem  \ref{l:timeTmarginal}.
\epf

\bpf[Proof of Theorem \ref{thm:uniquenessPractical}]
Theorem \ref{thm:uniqueness} would imply Theorem
\ref{thm:uniquenessPractical}, if $Q\{  \pi_0X \in\mathcal{A}_\infty(\rho,r)
\}=Q\{ \pi_+X \in \mathcal{B}_\infty(\rho,r) \} =1$
for every measure $Q$ which defines a stationary $\DV$-valued solution.
However this is precisely the content of
Theorems \ref{asymtotic behavior} and \ref{thm:avgFluc}.
\epf

\section{Basic Examples}
\label{sec:examples}

In this section we consider two simple examples to illustrate our
theory.  The first example is quite uniform in its behavior.  The
second one, though it is similar to the first, has a twist which makes
the estimates less uniform. It can be seen as a warm up for applying
our results to stochastically forced partial differential equations.
Finally in the next section we illustrate the point by considering
stochastically forced dissipative partial differential equations.

\subsection{A Uniform Example}

Consider the equation \eqref{sdde} with $X(t) \in \RR$ and
 $a(x)=-x(0)(1+\Psi(x))$ where
 \begin{equation*}
\Psi(x)=\int_{-\infty}^0 e^{-s^2+s} x(s)^2 ds  \ .
 \end{equation*}
 We now show that Theorem \ref{existence through Lyapunov} and Theorem
 \ref{thm:uniquenessPractical} apply, thus the system has a
 unique stationary solution.

First note that $V(x)=x(0)^2 + \Psi(x)^2$ is a Lyapunov function for
the system with $\gamma=1$, $\delta=\frac12$, $C_1=1$ and
$C_2=1$. To see this, apply It\^o's formula to  $V$ giving
\begin{equation*}
  dV(\pi_tX)= h(\pi_t X)
  dt + 2X(t)dW(t)
\end{equation*}
where
\begin{equation*}
  h(x)=-2x(0)^2+1+\int_{\R_-}(2s-1)e^{-s^2+s}ds\le -V(x)+1
\end{equation*}
and notice that $|x(0)| \leq V(x)^\frac12$. Next we establish
continuity in $C_\rho$, for any $\rho >0$, of the functional $a$.
For $x,\tilde x \in C_\rho$, a straightforward calculation gives
\begin{align}\label{eq:ex1Contractive}
  |\Psi(x) - \Psi(\tilde x)| \leq & \int_{-\infty}^0 e^{-s^2+s} \big( |x(s)|
   + |\tilde x(s)| \big) |x(s)-\tilde x(s)| ds\\
   \leq & \big( \|x\|_\rho + \|\tilde x\|_\rho  \big) \|x-\tilde x\|_\rho
   \int_{-\infty}^0 e^{-s^2+s}(1+|s|^\rho)^2 ds \notag
\end{align}
Since $|a(x)| \leq C(1+V(x))$ for some positive $C$ and all $x$
the existence of a stationary
solution is implied by
 Theorem \ref{existence through  Lyapunov} which
applies with $\beta=1$ and $\nu$ taken to be the delta
measure concentrated at zero.

To show uniqueness of the stationary solution we consider $|a(\pi_t
(x\cat y) ) - a(\pi_t (\tilde x\cat y) )|$ for $t>0$ where $x,\tilde x
\in \mathcal{A}_n(\frac34,\frac34)$ and $y \in
\mathcal{B}_n(\frac34,\frac34)$ where $\mathcal{A}_n$ and
$\mathcal{B}_n$ are as defined in Theorem
\ref{thm:uniquenessPractical}. From \eqref{eq:ex1Contractive} we have
\begin{align*}
  |a\big(\pi_t (x \cat y) \big) - a\big(\pi_t (\tilde x \cat y) \big)|
  \leq& e^{-t} \int_{-\infty}^0 e^{-|s|} \big( |x(s)|
  + |\tilde x(s)| \big) |x(s)-\tilde x(s)| ds\\
  \leq  & e^{-t}  \int_{-\infty}^0  e^{-|s|} 4n^2(1+|s|^\frac34)^2  ds
  <C(n)e^{-t} \ .
\end{align*}
Taking $\mathcal{K}_n(t)=C^2(n)e^{-2t}$ we can apply Theorem
\ref{thm:uniquenessPractical} and hence the stationary $\DV$-valued
solution is unique. It is important to stress
that Theorem~\ref{thm:uniquenessPractical} gives only uniqueness of
$\DV$-valued solutions which is natural since the very dynamics is defined only
on the set where  $a(x)=-x(0)(1+\Psi(x))$ is finite.

\subsection{A Less Uniform Example}

We now modify the previous pedagogical example making it less
uniform. We set $a(x)=-x(0)[1+\Psi(x)] +\Psi(x)^2$ where now  $\Psi$ equals
$\hat\Psi(x)$ if $\hat\Psi(x)< \infty$ and zero otherwise. Here $\hat\Psi(x)$
is given by
\begin{align*}
  \hat\Psi(x)=\int_{-\infty}^0 \exp\Big(-2|s|-\int_{-|s|}^0 x(r)dr\Big) x(s)^2 ds  \ .
\end{align*}
We again take   $V(x)=x(0)^2 + \hat \Psi(x)^2$ as our Lyapunov function
which produces
\begin{align*}
dV(\pi_tX) = [-2|X(t)|^2 - 4 \Psi(\pi_tX)^2 +1 ] dt + 2 X(t) dW(t)
\end{align*}
for $t> s$ when $\hat\Psi(\pi_sX) < \infty$. Since $-2|x|^2 - 4
\Psi(x)^2 + 1 \leq 1 -2V(x)$, we take $C_1=1$, $C_2=2$, $\beta=1$,
$\gamma=1$, $\delta=\frac12$ and $\nu$ a delta measure concentrated at
zero. Next observe that $\frac1s \int^0_{-|s|} X(r) dr \leq (\frac1s
\int^0_{-|s|} X(r)^2 dr)^\frac12\leq (\frac1s \int^0_{-|s|}
V(X(r))dr)^\frac12 $. Since the time average of the Lyapunov function
is less than $\frac{C_1}{C_2}=\frac12$ (see Theorem \ref{thm:avgFluc}
and Lemma \ref{l:fixPointNice}) we conclude that if $\hat\Psi(X(t_1))
< \infty$ then $\hat\Psi(X(t)) < \infty$ for all $t$. Showing that
$a(\ccdot)$ is locally Lipschitz on $C_\rho^-\cap \Nice[r]$ reduces to
showing that $\hat\Phi$ is locally Lipschitz on $C_\rho^-\cap
\Nice[r]$, for concreteness we choose $\rho=r=\frac34$. Define
$\Gamma(s,t)= \exp\left( -2(t-s)+ \int_{s}^t (x(s_1)+\tilde x(s_1))
  ds_1 \right)$ and take $x,\tilde x \in \mathcal{A}_n$ where
$\mathcal{A}_n$ is as in Theorem \ref{thm:uniquenessPractical}. Then
direct calculation gives
\begin{multline*}
  |\hat\Psi(x) - \hat\Psi(\tilde x)| \leq \int_{-\infty}^0 \Gamma(s,0) |x(s)-\tilde
  x(s)|[C+V(x(s))+V(\tilde x(s) )] ds\\
  \qquad \leq  C\|x-\tilde x\|_\frac34  \int_{-\infty}^0 \exp\Biggl\{
  -2|s|\Biggl(1-\sqrt{\frac12+n\frac{1+|s|^\frac34}{|s|}}\Biggr)\Biggr\} n(1+|s|^\frac34)^2 ds \ .
\end{multline*}
Hence Theorem \ref{existence through Lyapunov} implies the existence
of a stationary solution.  An analogous calculation give for $t>0$, $x,\tilde
x \in \mathcal{A}_n$ and $y \in \mathcal{B}_n$
\begin{align*}
  |\hat\Psi(\pi_t(x \cat y)) - \hat\Psi(\pi_t(\tilde x \cat y))| &\leq |\hat\Psi(x) -
  \hat\Psi(\tilde x)|\Gamma(0,t) \\
   &\leq K_n(t)=2n
   \exp\Biggl\{
  -2|t|\Biggl(1-\Bigl(\frac12+n\frac{1+|t|^\frac34}{|t|}\Bigr)^{1/2}\Biggr)\Biggr\}
 \end{align*}
which through Theorem \ref{thm:uniquenessPractical} give uniqueness of
the stationary $\DV$-valued solution.

\section{Application to Stochastically Forced Dissipative Partial
  Differential Equations}

Consider the stochastic differential equation
\begin{align}\label{eq:PDE}
  du(t)&=F(u(t))dt + GdB(t)\\
  u(t_0)&=u_0 \notag
\end{align}
on a Hilbert space $\X$ with a norm $\Norm{\ccdot}$ and inner product
$\ip{\ccdot}{\ccdot}$. Here $F$ maps a set $\tilde\X\subset\X$ into
$\X$, $B$ is an $\X$-valued cylindrical Brownian motion on a
probability space $\Theta$ and $G$ is a Hilbert-Schmidt operator
mapping the domain of $B$ into $\X$. We assume that $B(t)$ exists for
all positive and negative times. We also assume without further
comment that \eqref{eq:PDE} has strong, global, pathwise unique
solutions. The goal of this section is to show that under certain
assumptions this Markovian system on a possibly infinite dimensional
phase space can be reduced to finite dimensional system with memory
which has the same asymptotic behavior as the original system. In the
stochastic setting, the reduction was proved in \cite{b:Mattingly98b}.
The pathwise contraction of the small scales embodied in
\eqref{eq:dissipative} were used in \cite{b:Mattingly98}.  These and
related ideas were used to prove ergodic results for the stochastically
forced Navier-Stokes equations in
\cite{b:EMattinglySinai00,b:BricmontKupiainenLefevere01,b:KuksinShirikyan00}.
In \cite{b:Mattingly02,b:BricmontKupiainenLefevere02}, exponential
mixing is proved. The first of these uses a coupling construction and
gives explicit estimates on the dependence of initial conditions. A
different point of view and extension of these ideas can be found in
\cite{b:KuksinShirikyan02} and \cite{b:Mattingly03Pre}. We emphasize,
that the precise ergodic results presented here are not new. They fit
into the frameworks in
\cite{b:Mattingly02,b:Hairer02,b:KuksinShirikyan02,b:Mattingly03Pre}
which give further information about the rates of convergence.  Our
intention is to make explicit and clarify the idea of reducing the
system to one with memory which provides useful intuition.

We make three basic assumptions:
\begin{enumerate}
\item The dynamics admits a Markovian Lyapunov function. Namely, in
  complete analogy to our previous definition, there
  exist a function $U:\X \to\R\cup \{\infty\}$ so that
  \begin{enumerate}
  \item $\tilde\X\subset\DU$ 
  \item $U(u)\ge C_0 \Norm{u}^l$   for some $C_0,l>0$ and all $u
    \in \tilde\X$.
  \item If $u$ is a solution of equation \eqref{sdde} for $t \geq s$ 
    and $U(u(s)) < \infty$ then 
    \begin{equation*}
      dU(u(t)) = h(u(t))dt + f(u(t)) d\widetilde B(t).
    \end{equation*}
     Here $h:\tilde\X\to\R$ is a
    function satisfying
    \begin{equation*}
      h(u)<C_1 - C_2 U(u)^\gamma
    \end{equation*}
    for some constants $C_1,C_2,\gamma>0$ and $f:\tilde\X\to\R$ is a function
    satisfying
    \begin{equation*}
      |f(u)| \le C_3U(u)^\delta,
    \end{equation*}
    for some $\delta \in\left[0, (1+\gamma)/2\right)$ and $C_3>0$.
    Finally $\widetilde B$ is a standard
    Wiener process adapted to the flow generated by $(dB)$.
  \end{enumerate}

\item There exist a decomposition of $\X=\LL \oplus \HH$ where $\LL$
  is $d$-dimensional with $d < \infty$ with the property if $\Pl$ and
  $\Ph$ are the respective projections on $\LL$ and $\HH$ then
  \begin{multline}
    \label{eq:dissipative}
    \ip{F(u)-F(\tilde u)}{\Ph(u-\tilde u)}\leq \Norm{\Ph(u-\tilde
  u)}^2(-c_1+c_2U(u)^\gamma) \\+ c_3  \Norm{\Pl(u-\tilde
  u)}^{p_1}(1+ U(u)^{p_2} +  U(\tilde u)^{p_2})
  \end{multline}
  for all $u,\tilde u \in \tilde\X$ and some $c_i,p_i \geq 0$ with $\frac{c_1}{c_2} >
  \frac{C_1}{C_2}$.
\item $\Ph G B(t)=0$. This amounts to only saying that the subspace
  $\LL$ is forced.
\item The Markovian dynamics \eqref{eq:PDE}  possesses an invariant
  probability measure $\mu$ such that 
   $\mu\{ u \in \tilde\X \} =1$.
\end{enumerate}

\remark The last condition requiring the existence of an invariant
measure can almost certainly be derived from the first two. This would
be interesting in that it would replace the normal compactness
arguments with a dissipative argument much as the uniqueness theory
removes most topological considerations by using the dissipative
structure.

\remark In fact many of the results follow without modification even
if $\Ph GB \not = 0$. It
would require the consideration of drifts of the form $a(\pi_tX,t)$ in
\eqref{sdde} and further assumptions to control the behavior in $t$.
Though it is instructive to explore this option, for the PDE
applications it is more natural to stay in a Markovian framework,
shifting to the memory point of view only as needed. See for example
\cite{b:Mattingly02,b:KuksinPiatnitskiShirikyan02,b:KuksinShirikyan02,b:Mattingly03Pre}

\remark The asymmetry in the estimate on the differences in the $F$ is
an artifact of our memory setting. If we stayed in a more Markov setting it
would not be needed. See \cite{b:Mattingly02,b:Mattingly03Pre} for more details.

We can extend the invariant
measure $\mu$ to a stationary measure $\mathcal{M}$ on
$C((-\infty,\infty),\X) \times \Theta$. See \cite{b:EMattinglySinai00}
equation (6) for more details. Then using precisely the same
calculations as
in Theorem \ref{thm:avgFluc} and \ref{asymtotic behavior}
with $\mathcal{F}U(u,t)= | \int_0^t U(\pi_s x)^\gamma ds |
  - \frac{C_1}{C_2}|t| $ we obtain
\begin{theorem}\label{thm:PDEGrowth}
  Under the assumptions 1--4
  for any $\rho>0$ and $\varkappa>\frac12$
\begin{equation*}
 \mathcal{M}\left\{\sup_{t  } \frac{|\mathcal{F}U(u,t)|}{1+|t|^\varkappa} <
  \infty\right\}=\mathcal{M}\left\{\sup_{t } \frac{U(u(t))}{1+|t|^\rho} <
  \infty\right\}=1 \ .
\end{equation*}
\end{theorem}

We now prove the critical lemma which allows us to remap this problem
to the setting of the first half of the paper. Let $\mathcal{M}\Pl^{-1}$
 be the projection onto paths in
$C((-\infty,\infty),\LL) \times\Theta$.
Similarly we split equation \eqref{eq:PDE} into equations for
$(h(t),\ell(t)) \in \LL \oplus \HH$, obtaining
\begin{align}
    \frac{dh(t)}{dt}&=\Ph F(h(t)+\ell(t))\label{eq:h}\\
    d\ell(t)&=\Pl F(h(t)+\ell(t))dt + GdB(t)
\end{align}
We now show that there exists a function $\Psi:C^- \rightarrow
\HH$ so that the equation
\begin{align}\label{eq:reduced}
  d\ell(t) = \Pl F(\ell(t)+\Psi(\pi_t \ell))dt + \Pi_\ell G dB(t)
\end{align}
has the same asymptotic behavior as \eqref{eq:PDE}. This equation is
similar to what is called an inertial form in the theory of inertial
manifolds; there however, the function $\Psi$ depends only on the
present and not on the past. In some settings but not in all cases,
one can construct a stochastic inertial manifold (see
\cite{b:DaPratoDebussche96,b:ChueshovGirya94,b:BensoussanFlandoli95}).
However the constructions of this section work in most dissipative
settings. The reduced memory formulation \eqref{eq:reduced} is akin to
the reduction of a dynamical system done in the context of symbolic
dynamics. By only having some coarse description of the dynamics, but
for a time interval of infinite length, one can reconstruct the exact
position. Usually the symbolic dynamics encodes the forward
trajectory. Here we are encoding the state of some subset of the
variables in the infinite past.
\begin{theorem}\label{thm:mainReduction}
  Under the same assumptions
  there exists an $\HH$-valued function $\Psi$ defined on $C((-\infty,0],\LL)$  so
  that the following holds:
 \begin{enumerate}
 \item For $\mathcal{M}\Pl^{-1}$-almost every $(\ell_{(-\infty,0]},B)\in
   C((-\infty,0],\LL) \times \Theta$ if
   $\Psi(\pi_t\ell_{(-\infty,0]})=h(t)$ with $t \leq 0$ then
   $u(t)=(\ell(t),h(t))$ is a solution to \eqref{eq:PDE} with
   noise realization $B(t)$.
 \item If $\Psi_{s,t}(\ell_{(-\infty,0]},h_0)$ is the solution to
   \eqref{eq:h} at time $t$ with initial data $h_0 \in \HH$ at time
   $s$ and exogenous forcing $\ell$, then for
   $\mathcal{M}\Pl^{-1}$-almost every $\ell_{(-\infty,0]}$ one has
   $\lim_{s\phi - \infty} \Psi_{s,t}(\ell_{(-\infty,0]},h_0)=
   \Psi(\pi_t \ell_{(-\infty,0]})$ for any $h_0$.
 \item Fix $r>0$ and $\kappa\in(\frac12,1)$. There exists a constant
   $C(r,\kappa,n)$ so that for $\mathcal{M}\Pl^{-1}$ almost
   every $\ell, \tilde \ell \in A_n(\kappa)$,
   $\Norm{\Psi(\pi_0\ell)-\Psi(\pi_0\tilde \ell)}^2\leq C(r,\kappa,n)\|\ell -
   \tilde \ell\|_r^{p_1}$ where
   \begin{align*}
     A_n(\kappa)=\left\{ \ell \in C\big((-\infty,0],\LL\big) : \sup_{s \leq 0} \frac{
     U(u(s)) + \mathcal{F}U(u,s)}{1 + |s|^\kappa} < n \right\}
   \end{align*}
   and $u(s)=\ell(s)+ \Psi(\pi_s \ell)$. Furthermore there exist
   $\gamma_n >0$ and $K_n>0 $so that if $t>0$ and $\ell(s)=\tilde
   \ell(s)$ for $s \in [0,t]$ then
   $\Norm{\Psi(\pi_t\ell)-\Psi(\pi_t\tilde \ell)}^2\leq K_n\exp(-\gamma_n t)$
 \end{enumerate}
\end{theorem}

\bpf For $\mathcal{M}\Pl^{-1}$-almost every $(\ell(t),B(t))$ there is a
corresponding $h(t)$ so that $u(t)=(\ell(t),h(t))$ is a solution to
\eqref{eq:PDE} with forcing $B(t)$. For $s\leq 0$, $\tilde
h(s)=\Psi_{-t,s}(\ell_{(-\infty,0]},\tilde h_0)$. Defining
$\rho(s)=h(s)-\tilde h(s)$ we have
\begin{align*}
  \frac{d \Norm{\rho(s)}^2}{dt} &=\ip{F(\ell(s)+h(s)) -
  F(\tilde\ell(s)+h(s))}{\rho(s)}\\
&\leq  \Norm{\rho(s)}^2(-c_1+c_2U(u(s))^\gamma) \ .
\end{align*}
Hence using Theorem \ref{thm:PDEGrowth} to continue
\begin{align*}
  \Norm{\rho(s)}^2 &\leq ( \Norm{\tilde h_0} + \Norm{h(-t)})^2
  \exp\left( -c_1 (|t| -|s|) + c_2\int^s_{-t} U(u(\tau))^l d\tau
  \right)\\
  & \leq ( \Norm{\tilde h_0} + \|h\|_\kappa(1+|t|^\kappa )^2 \exp\left(
     -(c_1 - c_2\frac{C_1}{C_2}) (|t| -|s|)+ C(\kappa)(1+|t|^\kappa)  \right)
\end{align*}
for some $\kappa \in (\frac12,1)$. Since by assumption $c_1 -
c_2\frac{C_1}{C_2} > 0$, $\Norm{\rho(s)} \rightarrow 0$ as $t
\rightarrow - \infty$. This proves the first and second claim of the
theorem. To see the third consider two pairs of solution $(\ell,B)$
and $(\tilde \ell,\tilde B)$ in $A_n$.  For $ \mathcal{M}\Pl^{-1}$-almost every such
pair, in the same way as the previous estimate one obtains
\begin{multline*}
  \Norm{\rho(0)}^2 \leq c_3\|\ell-\tilde \ell\|_r^{p_1} \int^0_{-\infty} \exp\left( -(c_1 - c_2
  \frac{C_1}{C_2})|t| + nc_2(1+|t|^\kappa)\right)\times \\(1+|t|^r)^{p_1}
  (1+2n^{p_2}(1+|t|^\kappa)^{p_2})dt.
\end{multline*}
\epf

With Theorem \ref{thm:mainReduction} in hand, we can define dynamics
with memory on finite dimensional space $\LL$ which is isomorphic to
$\RR^d$. We set $Q=\mathcal{M}\Pl^{-1}$ and consider an equation of the
form \eqref{sdde} with $a(x)=\Pl F(x(0)+\Phi(x))$ and with a Lyapunov
function $V(x)=U(x(0)+\Psi(x))$. We then arrive at an equation of the
form \eqref{sdde}. Hence by invoking theorem
\ref{thm:uniquenessPractical}, we obtain the following result.
\begin{theorem}\label{thm:PDEErg}
  Assuming all of the assumptions of this section and additionally
  that $\Pl G$ is of full rank and that
  \begin{equation}\label{eq:lowEstimate}
    |\Pl F(\ell +  h)-\Pl F(\ell + \tilde h)|^2\leq c_4( 1 + U(\ell +
      h)^{p_3}+ U(\ell +  \tilde h)^{p_3})\Norm{h-\tilde
      h}^{p_4} \ .
  \end{equation}
  for all $\ell \in \LL$ and $h,\tilde h \in \HH$ such that 
  $\ell+h,\ell+\tilde h\in \tilde\X$ and some $c_4 \geq
  0$, and $p_i \geq 0$. Then the invariant solution $\mathcal{M}$ is
  the unique $\DU$-valued one.
\end{theorem}

Notice that ergodicity only requires that all of the modes up to a
given scale are forced. In other words, it is sufficient for the
system to be elliptic only up to a certain scale to ensure ergodicity.
The arguments used are ``soft'' in that they do not require explicit
geometric information as hypoelliptic arguments do. For this reason,
it is reasonable to call the system ``effectively elliptic'' because
the reduced memory equation is truly elliptic and hence the arguments
are relatively ``soft.''

\bpf Notice that the sets in Theorem \ref{thm:uniquenessPractical} are
just the projects of the $t\leq 0$ and $t\geq 0$ parts of the set
defined in Theorem \ref{thm:PDEGrowth} and agree with the $A_n$
defined above. Let $x_1$, $x_2$, $y$ be chosen as in Theorem
\ref{thm:uniqueness}. Since $a(x)=\Pl F(x(0)+ \Psi(x))$,
\begin{align*}
  |a(&\pi_t(x_1\cat y))-a(\pi_t(x_2 \cat y))|^2\\&= |\Pl F(y(t)+
  \Psi(\pi_t(x_1 \cat y)))-\Pl F(y(t)+ \Psi(\pi_t(x_2 \cat y)))|^2\\
  & \leq c_4( 1 + V(\pi_t(x_1 \cat y))^{p_3}+ V(\pi_t(x_2\cat y))^{p_3})\Norm{
    \Psi(\pi_t(x_1 \cat y))- \Psi(\pi_t(x_2 \cat y))}^{p_4} \\
  & \leq c_4(1+ 2n^{p_3}(1+|t|^\kappa)^{p_3})K_n\exp(-p_4\gamma_n t)
\end{align*}
As this bound is integrable, Theorem \ref{thm:uniquenessPractical}
completes the proof.
\epf

One interesting consequence of the memory point of view is following
factorization of  the invariant measure into a measure living on the path
space $C((-\infty,0],\LL)$ and an atomic measure living in $\HH$ which
depends only on the choice in  $C((-\infty,0],\LL)$. This
factorization show that the random attractor projected into  $\HH$
space is a single point attractor fibered over the choice of
trajectory in $C((-\infty,0],\LL)$.
\begin{theorem}\label{thm:factor}
  Assuming all of the conditions of Theorem \ref{thm:PDEErg}
  hold. Then the following factorization of the invariant
  measure $\mu$ holds: for any $A \subset \X$
  \begin{align*}
    \mu(A) = \int \ONE_A(\ell(0)+\Psi(\ell))
    \mathcal{M}\Pl^{-1}(d\ell) \ .
  \end{align*}
\end{theorem}

\subsection{The 2D Stochastic Navier Stokes Equation}

Consider the incompressible Navier Stokes equation with mean zero flow
on the two dimensional unit torus,$\T^2$, agitated by a stochastic forcing with no
mean flow. By projecting out the pressure, we obtain the following
It\^o equation for $u(x,t)=(u^{(1)}(x,t),u^{(2)}(x,t)) \in \X=L^2(\T^2)
\times L^2(\T^2)$
\begin{align*}
  du(x,t)=\left[\nu \Delta u + B(u,u)\right] dt + G dW(t)
\end{align*}
where $B(u,v)=P_{div}(u\cdot \nabla ) u$, $P_{div}$ is the projection
onto divergence-free vector fields, $G$ a Hilbert-Schmidt operator
mapping the cylindrical Weiner process $W$ into $\X$. We assume that
there exists $\mathcal{K}_{\cos}$, $\mathcal{K}_{\sin} \subset \Z^2$
so that $Im(G)=\LL=\mbox{span}( \sin(2\pi k\cdot x),\cos(2\pi m\cdot
x) : k \in \mathcal{K}_{\sin}, m \in \mathcal{K}_{\cos})$. We define
$N_0$ to be the largest integer multiple of $2\pi$ so that if
$2\pi|k| < N_0$ then $k \in \mathcal{K}_{\cos} \cap
\mathcal{K}_{\sin}$. Similarly we define $N_1$ to be the smallest
integer so that if $2\pi|k| > N_1$ then  $k \not \in \mathcal{K}_{\cos} \cup
\mathcal{K}_{\sin}$. We want to show that if $N_0$ is sufficiently
large then Theorem \ref{thm:PDEErg} and \ref{thm:factor} hold. We take
$U(u)=\Norm{ \nabla u}^2$ as the Lyapunov function. Standard
calculations  show that $U$  satisfies our conditions for a Lyapunov function with
$C_1=\text{Tr}GG^*$, $C_2=2\nu$, and $C_3=2$.  (See the enstrophy calculations
in \cite{b:EMattinglySinai00}). Lastly, we need to verify
\eqref{eq:dissipative} and \eqref{eq:lowEstimate}. We begin with the
first, setting $F(u)=\nu \Delta u + B(u,u)$, $\HH=\X / \LL$ and $\rho=
u-\tilde u$ one has
\begin{align*}
  \ip{F(u)-F(\tilde u)}{\Ph (u-\tilde u)}&\leq -\nu \Norm{\Ph
    \nabla(u-\tilde u)}^2 + \ip{B(\rho,u)}{\Ph (u-\tilde u)} \\&\qquad
  + \ip{B(u,\rho)}{\Ph
    (u-\tilde u)}\\
  & \leq \left(-\frac{\nu N_0^2}{2} + \frac{C}{\nu}\Norm{\nabla u}^2
  \right)\Norm{\Ph (u-\tilde u)}^2\\ &\qquad  + \frac{CN_0^2}{\nu}\left(\Norm{\nabla
    u}^2+ \Norm{\nabla \tilde u}^2\right) \Norm{\Pl (u-\tilde u)}^2
\end{align*}
where the constant $C$ is independent of $N_0$. Hence if $\frac{\nu
  N_0^2}{2} > \frac{C}{2\nu^2}\text{Tr}GG^*$, the assumption in
  \eqref{eq:dissipative} holds.

  Lastly we check the condition used to control the paths in $\LL$.
  For any $h,\tilde h \in \HH$, $\ell \in \LL$, and $v \in \LL$ with
  $\Norm{v}=1$ one has
\begin{align*}
  \ip{F(\ell + h) - F(\ell + \tilde h)}{v} \leq & \ip{B(h-\tilde
    h,\ell + h)}{v} + \ip{B(\ell + \tilde h,h-\tilde h)}{v} \\
  \leq& C \Norm{h- \tilde h} \Norm{ \nabla^2 v} \left( \Norm{ \nabla
      (\ell + h)} +\Norm{ \nabla
      (\ell + \tilde h)} \right) \\
  \leq & C N^2_1 \Norm{h- \tilde h}\left( \Norm{ \nabla (\ell + h)}
    +\Norm{ \nabla (\ell + \tilde h)} \right)
\end{align*}
Hence \eqref{eq:lowEstimate}, holds and we obtain the following
result.
\begin{theorem}
  Let $C$ be the constant from the above calculations. If $N_0^2 >
 C \frac{\text{Tr}GG^*}{\nu^3}$ then Theorems  \ref{thm:PDEErg} and
  \ref{thm:factor} apply to the Stochastic Navier-Stokes equation.
\end{theorem}
\subsection{The 1D Stochastic Ginsburg-Landau Equation}

As a second PDE example, we consider the stochastically forced
Cahn-Allen/Ginsburg-Landau equation in a one dimensional periodic
domain. Ergodic results for this equation have been proved in
\cite{b:ELui02,b:Hairer02}. The general framework of
\cite{b:Mattingly02} applies equally well to this setting.

Consider
\begin{align}
  \label{eq:GL}
  du(x,t)=\left[\nu \Delta u + u - u^3\right] dt +dW(x,t)
\end{align}
where $W(x,t)=\sum_{\mathcal{K}} e_k(x) \sigma_k \beta(t)$, $\beta_k$
are independent standard Brownian motions, $\sigma_k$ are positive
constants and $e_k$ are the elements of the real Fourier basis $\{ 1,
\sin(2\pi x), \cos(2\pi x), \sin(4\pi x) ,\cos(4\pi x), \cdots \}$. We
denote by $\lambda_k$ the eigenvalue of $-\Delta$ associated with
$e_k$. As before we consider the case when $\sigma_k>0$ for only a
finite number of $k$ and define $N_0$ by the smallest integer so that
if $\lambda_k< 4\pi^2 N_0^2$ then $\sigma_k >0$.

As in the previous example, we let $\LL$ be the span of the $e_k$ with
$\sigma_k>0$ and $\HH$ the span of the remaining $e_k$. Then
$\X=L^2([0,1])=\LL \oplus \HH$ and by $\Norm{\cdot}$ we mean the
$L^2$-norm on $\X$.

We use the Lyapunov function $U(u)=\Norm{u}^2 + \Norm{\nabla u}^2$.
Direct calculation gives (see for instance \cite{b:ELui02,b:Sm94})
\begin{multline*}
  dU(t) = 2\left[ \Norm{u}^2 - \Norm{\Delta u}^2 + (1 -\nu)\Norm{\nabla
    u}^2 - 3\nu \Norm{u \nabla u}^2 - \int u^4 dx +
  \frac12\mathcal{E}_0 + \frac12\mathcal{E}_1 \right] dt \\+ 2\ip{u}{dW} -
  2\ip{\Delta u}{dW} .
\end{multline*}
where $\mathcal{E}_m = \sum \lambda_k^m \sigma_k^2$.  Using Jensen's
inequality on the $L^4$ norm, the fact that $\frac12K_0 - x^2 > x^2
-x^4 $ for some positive $K_0$, $\Norm{\Delta u}^2 >4 \pi^2
\Norm{\nabla u}^2$ and the non-positivity of the third term, we see
that $U$ is a Lyapunov function with $C_1=K_0 + \mathcal{E}_0
+\mathcal{E}_1$, $C_2=2$ , $\gamma=1$, $\delta=\frac12$ and
$C_3=2\sigma_{max}^2(0) + 2\sigma_{max}^2(1)$. Here $\sigma_{max}^2(m)
= \max \sigma_k^2 \lambda_k^m$.

All that remains is to prove the estimates \eqref{eq:dissipative} and
\eqref{eq:lowEstimate} hold. To see the first recall that the Sobolev
embedding theorem implies that $|f|_{L^\infty}^2 < CU(f)$ and hence if
$F(u)=\Delta u + R(u)$ and $R(u)=u-u^3$ then
\begin{align*}
  \ip{F(u)-F(\tilde u)}{\Pi_h (u - \tilde u)} &= -\Norm{\nabla \Pi_h
    (u - \tilde u)}^2 +\ip{R(u)-R(\tilde u)}{u - \tilde u}
  \\ &\qquad + \ip{R(u)-R(\tilde u)}{\Pi_\ell (u - \tilde u)} \\
  & \leq (1-4\nu \pi^2 N_0^2) \Norm{\Pi_h (u - \tilde u)}^2
  +\Norm{\Pi_\ell (u - \tilde u)}^2 \\&\qquad + |\ip{R(u)-R(\tilde
    u)}{\Pi_\ell (u - \tilde u)}| \ .
\end{align*}
To finish this estimate, observe that
\begin{align*}
  |\ip{R(u)-R(\tilde
      u)}{\Pi_\ell (u - \tilde u)}| & \leq ( 1+ 2 |u|^2_{L^\infty} +2
      |\tilde u|^2_{L^\infty}) \Norm{\Pi_\ell (u - \tilde u)}^2\\
      & \leq C( 1 + U(u) + U(\tilde u) )  \Norm{\Pi_\ell (u -
      \tilde u)}^2 \ .
\end{align*}
To see \eqref{eq:lowEstimate}, set $u=\ell + h$ and $\tilde u = \ell +
\tilde h$ and observe that
\begin{align*}
  \Norm{\Pi_\ell(F(\ell + h)-F(\ell + \tilde h)) } = & \Norm{(\ell +
    h)^3
    - (\ell + \tilde h)^3}+ \Norm{h-\tilde h}\\
    = & \Norm{\int_u^{\tilde u} 3 v^2 dv }+\Norm{h-\tilde h}\\
    \leq& \Norm{(3u^2 +3\tilde u^2)|h-\tilde h|}+\Norm{h-\tilde h} \\
  \leq & C(|u|^2_{L^\infty}+|\tilde u|^2_{L^\infty}+1) \Norm{h - \tilde h} \\
  \leq & C \left(U(u) + U(\tilde u)+1 \right) \Norm{h - \tilde h}.
\end{align*}
In the last line we have used the
Sobolev embedding theorem.

In light of the above calculations, we have
\begin{theorem}
  If $N_0^2 > \frac{1}{4 \pi^2 \nu}$
  then Theorems \ref{thm:PDEErg} and \ref{thm:factor} apply to the
  stochastically forced  Ginsburg-Landau Equation \eqref{eq:GL}.
\end{theorem}

\section{Acknowledgments}
We would like to thank Yakov Sinai, Toufic Suidan and Alexander
Veretennikov for their input and comments. We would also like to thank
the hospitality and support of the Institute for Advanced Study in
Princeton where this work was undertaken. The second author was
partial supported by the NSF under grant DMS-9971087. And both authors
received support from the NSF under grant DMS-9729992.

\appendix

\section{Comparison of Measures on Path Space}
\label{sec:measurePath}
Suppose that we have two measures $P^{(1)}$ and $P^{(2)}$ on the space
$C^+\times C^+$ which define solutions for equations
\begin{equation}
\begin{split}
dX_i(t)&=f_i(t,X_i[0,t])dt+dW(t),\ t\ge 0,\quad i=1,2\ \mbox{\rm respectively},\\
X_i(0)&=x_0.
\end{split}
\label{2SDDEs}
\end{equation}

Here for fixed $t$ functions
$f_1$ and $f_2$ map the space $C_{[0,t]}=C([0,t],\RR^d)$ to $\RR^d$.
By $X[0,t]$ we mean the segment of the trajectory on $[0,t]$.
Let $T\in (0,\infty]$ and $\mathcal{B}\subset C_{[0,T]}$.
Define measures $P^{(i)}_{[0,T]}(\ccdot ; \mathcal{B})$ on the path space as:
\begin{equation*}
  P^{(i)}_{[0,T]}(A; \mathcal{B}) =
  P\{X_i[0,T]\in A\cap\mathcal{B}\},\ \mbox{\rm for}\ A \subset C_{[0,T]}.
\end{equation*}
Also define $D(t,\ccdot)=f_1(t,\ccdot)-f_2(t,\ccdot)$.

In this
setting, we have the following result which is a variation on Lemma
B.1 from \cite{b:Mattingly02}.
\begin{lemma} \label{l:compareMeasures}
  Assume there exists a constant $D_*\in(0,\infty)$ such that
  \begin{align}\label{eq:Novikov}
     \exp\left\{\frac12 \int_0^T \big|D\big(t,X[0,t]\big)\big|^2 dt \right\}
    \ONE_\mathcal{B}(X[0,t])  < D_*
  \end{align}
   almost surely with respect to both measures $P^{(1)}$ and $P^{(2)}$.
  Then the measures
  $P^{(1)}_{[0,T]}(\ccdot ; \mathcal{B})$ and $P^{(2)}_{[0,T]}(\ccdot ; \mathcal{B})$ are
  equivalent.

  \end{lemma}

  \bpf
  Define the auxiliary SDEs
  \begin{align*}
    dY_i(t) &=f_i\big(t, Y_i[0,t]\big) \ONE_{\mathcal{B}(t)}(Y_i[0,t]) dt + dW(t)
  \end{align*}
  where
  $
    \mathcal{B}(t)=\{ x \in C_{[0,t]}: \exists \bar{x} \in
    \mathcal{B} \mbox{ such that } x(s)=\bar{x}(s) \mbox{ for $s
    \in[0,t]$}  \}
  $.
  Solutions $Y_i(t)$ to these equation can be constructed as
  \begin{equation*}
    Y_i(t)=X_i(t)\ONE_{\{t\le\tau\}}+[W(t)-W(\tau)+X_i(\tau)]\ONE_{\{t>\tau\}}.
  \end{equation*}
  Here $(X_i(t),W)$ is the solution to equation (\ref{2SDDEs}) and
  $\tau=\inf\{s>0: X_i[0,s]\not\in \mathcal{B}(s)\}$.

  Denote
  $D_\mathcal{B}(t,x) =[f_1(t,x) - f_2(t,x)]\ONE_{\mathcal{B}(t)}(x)$.
  The assumption on $D$ in \eqref{eq:Novikov} and the definition of
  $\mathcal{B}(t)$ imply that
  \begin{equation*}
    \exp\left\{\frac12\int_0^T \big| D_\mathcal{B}
    \big(t,X[0,t]\big)\big|^2 dt
    \right\}<D_*\quad \mbox{a.s.}
  \end{equation*}
  under both measures $P^{(i)}_{Y[0,t]}$ defining
  solutions to auxiliary equation with $i=1$
  and $i=2$.
  Hence Novikov's condition is satisfied for the difference of the drifts
  of the auxiliary equations and
  the Girsanov theorem implies that
  $\frac{dP^{(1)}_{Y[0,t]}}{dP^{(2)}_{Y[0,t]}}(x) = \mathcal{E}(x)$ where
  the Radon--Nikodym derivative evaluated at a trajectory $x$ is defined by
  the stochastic exponent:
  \begin{align*}
    \mathcal{E}(x)=\exp\left\{\int_0^T
     \left\langle D_{\mathcal{B}}(s, x[0,s]), dW(s)\right\rangle - \frac12
      \int_0^T |D_{\mathcal{B}}(s, x[0,s])|^2 ds \right\}.
  \end{align*}

  Note that restrictions of the measures $P^{(i)}_{Y[0,t]}$ on the set
  $\mathcal{B}$ coincide with $P^{(i)}_{[0,t]}(\ccdot,\mathcal{B})$.
  This proves that $P^{(1)}_{[0,t]}(\ccdot,\mathcal{B})$ is absolutely
  continuous with respect to $P^{(2)}_{[0,t]}(\ccdot,\mathcal{B})$.
  The reverse relation follows by symmetry and the proof is complete.
  \epf

\section{Cauchy problem}\label{app:Cauchy problem}
In this appendix we study existence and
uniqueness of solution to the Cauchy problem for \eqref{sdde}.
\begin{theorem}\label{thm:Cauchy existence}
Suppose the functional $a(\ccdot)$ is locally Lipschitz on $C_\rho$
with respect
to the norm $\|\ccdot\|_{\rho}$. If
$x\in C_{\rho}$ then for any realization of standard
Wiener process $W$ on any probability space
there exist a positive stopping time $T$ and
a continuous process $X(t)$ on the same probability space
with the following properties:
\begin{enumerate}
\item $X(t)=x(t)$ for $t\le0$ almost surely.
\item The couple $(X,W)$ solves equation \eqref{sdde} on $[0,T]$.
\item The process $X$ is adapted to the flow generated by $W$.
\end{enumerate}
Any other process with this properties coincides with $X$ almost surely.
\end{theorem}
\bpf
Fix a trajectory of $W$ on $[0,\infty)$ and for any positive
$T$ define operator $\Phi:C_{[0,T]}\to C_{[0,T]}$ by
\begin{equation*}
\Phi(y)(t)=x(0)+\int_0^ta(\pi_s(x \cat y))ds+W(t).
\end{equation*}

Then
\begin{equation*}
  |\Phi(y_1)(t)-\Phi(y_2)(t)|\le
  \int_0^t|a(\pi_s(x \cat y_1))-a(\pi_s(x\cat y_2))|ds.
\end{equation*}
If $\sup_{s\in[0,T]}|y_i(s)|\le M, i=1,2$ for some $M>0$,
then
\begin{equation*}
 \|\pi_s(x \cat y_i)\|_{\rho}\le M + \|x\|_\rho,\quad i=1,2.
\end{equation*}

Hence
\begin{equation*}
 |a(\pi_r(x \cat y_1))-a(\pi_r(x \cat y_2))|<K
 \|\pi_r(x\cat y_1)-\pi_r(x \cat y_2)\|_\rho
\end{equation*}
where $K=K(M+\|x\|_\rho)$ is the local Lipschitz constant.
If $|y|$ is bounded by $M$
then
\begin{equation*}
 \|\pi_r(x \cat y_1)-\pi_r(x \cat y_2)\|_\rho\le\sup_{s\in[0,t]}|y_1(s)-y_2(s)|(1+t^{\rho})
\end{equation*}
implies that
\begin{equation*}
  |\Phi(y_1)(t)-\Phi(y_2)(t)|\le t\cdot K(1+t^{\rho})
  \sup_{s\in[0,t]}|y_1(s)-y_2(s)|
  \quad\mbox{for\ } t\in[0,T].
\end{equation*}
Choose $M>|x_0|+2$ and
\begin{equation*}
T=\sup\Big\{t:tK(1+t^{\rho})<\frac12,
|W(t)|<1\mbox{\ and\ } Ct<1 \Big\}\wedge t_0
\end{equation*}
where $C=C(M)>0$ and $t_0>0$  are such that for $t<t_0$
the drift $|a(\pi_t(x \cat y))|$ is bounded by $C$
if $|y|$ is bounded by $M$.

Then $\Phi$ is a contraction on the $L^\infty$ maps from compact set
\begin{multline*}
  D=\biggl\{y\in C_{[0,T]}:  y(0)=0,\sup_{s\in[0,T]}|y(s)|\le
  M,\\ w_\delta(y)\le C\delta+\omega_\delta(W)\ \mbox{for all\ }
  \delta>0\biggr\}
\end{multline*}
to itself where $w_\delta(y)$ is the $\delta$-modulus of continuity of $y$.

There exists unique fixed point in $D$ of this map which gives the desired solution.
Since the choice of $M$ is arbitrary, we conclude that this solution is unique.
\epf

When it applies, Theorem \ref{thm:Cauchy existence}  implies that the
dynamics corresponding to the equation \eqref{sdde}
is defined at least up to some random moment $T$.

Let $T_\infty$ be the largest time such that a solution exists on
$[0,T_\infty)$.
The question of global existence will be answered if one shows that
if $T_\infty$ is finite then the solution can be in fact
extended beyond it.

\begin{theorem}
If there exists a Lyapunov function (see Section \ref{sec:Definitions}) for the dynamics corresponding
to \eqref{sdde} built in Theorem \ref{thm:Cauchy existence} then
the time $T$ in this theorem can be chosen to be equal to $\infty$,
i.e. pathwise uniqueness and strong existence hold globally.
\end{theorem}
\bpf
Essentially we need to show that $X(t)$ does not escape to infinity
in finite time a.s.
Introduce the stopping time $\tau_R=\inf\{t>0: V_t\ge R\}$.
Then by the definition of Lyapunov function we have
\begin{equation*}
V(\tau_R\wedge t)\le Ct+\int_0^{\tau_R\wedge t}h(\pi_s(X))dW(s).
\end{equation*}
Hence
\begin{equation*}
  \E V(\tau_R\wedge t)\le Ct.
\end{equation*}
This inequality with $\E V(\tau_R\wedge t) > R\cdot \PP\{\tau_R\le t\} $
implies that for any  $t>0$
\begin{equation*}
  \PP\{\mbox{for every\ } R>0\ \mbox{there exists \ } s\le t\ \mbox{such that\ } |V_s|>R\}=0. \end{equation*}
So $V_s$ is finite for all $s\le t$
which implies that $X(s)$ is finite for all $s\le t$.
\epf


\begin{thebibliography}{DPD96}
\bibitem[Bak02]{b:Bakhtin02}
Yu.~Yu. Bakhtin.
\newblock Existence and uniqueness of stationary solution of %nonlinear
  stochastic differential equation with memory.
\newblock {\em Theory Probab. Appl}, 47(4):764--769, 2002.

\bibitem[BF95]{b:BensoussanFlandoli95}
Alain Bensoussan and Franco Flandoli.
\newblock Stochastic inertial manifold.
\newblock {\em Stochastics Stochastics Rep.}, 53(1-2):13--39, 1995.

\bibitem[Bil68]{b:Billingsley68}
Patrick Billingsley.
\newblock {\em Convergence of probability measures}.
\newblock John Wiley \& Sons, Inc., New York-London-Sydney, 1968.

\bibitem[BKL01]{b:BricmontKupiainenLefevere01}
J.~Bricmont, A.~Kupiainen, and R.~Lefevere.
\newblock Ergodicity of the 2{D} {N}avier-{S}tokes equations with random
  forcing.
\newblock {\em Comm. Math. Phys.}, 224(1):65--81, 2001.
\newblock Dedicated to Joel L. Lebowitz.

\bibitem[BKL02]{b:BricmontKupiainenLefevere02}
J.~Bricmont, A.~Kupiainen, and R.~Lefevere.
\newblock Exponential mixing of the 2{D} stochastic {N}avier-{S}tokes dynamics.
\newblock {\em Comm. Math. Phys.}, 230(1):87--132, 2002.

\bibitem[CG94]{b:ChueshovGirya94}
I.~D. Chueshov and T.~V. Girya.
\newblock Inertial manifolds for stochastic dissipative dynamical systems.
\newblock {\em Dopov./Dokl. Akad. Nauk Ukra\"\i ni}, 7:42--45, 1994.

\bibitem[DPD96]{b:DaPratoDebussche96}
Giuseppe Da~Prato and Arnaud Debussche.
\newblock Construction of stochastic inertial manifolds using backward
  integration.
\newblock {\em Stochastics Stochastics Rep.}, 59(3-4):305--324, 1996.

\bibitem[EL02]{b:ELui02}
Weinan E and Di~Lui.
\newblock Gibbsian dynamics and invariant measures for stochastic dissipative
  {PDE}s.
\newblock {\em Journal of Statistical Physics}, 108(5/6):1125--1156, 2002.

\bibitem[EMS01]{b:EMattinglySinai00}
Weinan E, J.~C. Mattingly, and Ya~G. Sinai.
\newblock Gibbsian dynamics and ergodicity for the stochastic forced
  navier-stokes equation.
\newblock {\em Comm. Math. Phys.}, 224(1), 83--106, 2001.

\bibitem[FP67]{b:FoiasProdi67}
C.~Foia{\c{s}} and G.~Prodi.
\newblock Sur le comportement global des solutions non-stationnaires des
  \'equations de {N}avier-{S}tokes en dimension {$2$}.
\newblock {\em Rend. Sem. Mat. Univ. Padova}, 39:1--34, 1967.

\bibitem[Hai02]{b:Hairer02}
Martin Hairer.
\newblock Exponential mixing properties of stochastic pdes through asymptotic
  coupling.
\newblock {\em Probab Theory Relat Fields}, 124(3):345--380, 2002.

\bibitem[IN64]{b:ItoNisio}
Kiyoshi It\^o and Makiko Nisio.
\newblock On stationary solutions of a stochastic differential equation.
\newblock {\em J. Math. Kyoto Univ.}, 4:1--75, 1964.

\bibitem[KPS02]{b:KuksinPiatnitskiShirikyan02}
Sergei Kuksin, Andrey Piatnitski, and Armen Shirikyan.
\newblock A coupling approach to randomly forced nonlinear {PDE}s. {II}.
\newblock {\em Comm. Math. Phys.}, 230(1):81--85, 2002.

\bibitem[KS00]{b:KuksinShirikyan00}
Sergei Kuksin and Armen Shirikyan.
\newblock Stochastic dissipative {P}{D}{E}s and {G}ibbs measures.
\newblock {\em Comm. Math. Phys.}, 213(2):291--330, 2000.

\bibitem[KS02]{b:KuksinShirikyan02}
Sergei Kuksin and Armen Shirikyan.
\newblock Coupling approach to white-forced nonlinear {PDE}s.
\newblock {\em J. Math. Pures Appl. (9)}, 81(6):567--602, 2002.

\bibitem[Mat98]{b:Mattingly98b}
Jonathan~C. Mattingly.
\newblock {\em The Stochastically forced {N}avier-{S}tokes equations: energy
  estimates and phase space contraction}.
\newblock PhD thesis, Princeton University, 1998.

\bibitem[Mat99]{b:Mattingly98}
Jonathan~C. Mattingly.
\newblock Ergodicity of $2${D} {N}avier-{S}tokes equations with random forcing
  and large viscosity.
\newblock {\em Comm. Math. Phys.}, 206(2):273--288, 1999.

\bibitem[Mat02]{b:Mattingly02}
Jonathan~C. Mattingly.
\newblock Exponential convergence for the stochastically forced
  {N}avier-{S}tokes equations and other partially dissipative dynamics.
\newblock {\em Comm. Math. Phys.}, 230(3):421--462, 2002.

\bibitem[Mat03]{b:Mattingly03Pre}
Jonathan~C. Mattingly.
\newblock On recent progress for the stochastic {N}avier {S}tokes equations.
\newblock In {\em Journ\'ees \'Equations aux d\'eriv\'ees partielles},
  Forges-les-Eaux, XI:1--52, 2003.
\newblock see
  http://www.math.sciences.univ-nantes.fr/edpa/2003/html/.

\bibitem[MY02]{b:MasmoudiYoung02}
Nader Masmoudi and Lai-Sang Young.
\newblock Ergodic theory of infinite dimensional systems with applications to
  dissipative parabolic {PDE}s.
\newblock {\em Comm. Math. Phys.}, 227(3):461--481, 2002.

\bibitem[Pro90]{b:Protter90}
Philip Protter.
\newblock {\em Stochastic Integration and Differential Equations: a new
  approach}.
\newblock Springer-Verlag, 1990.

\bibitem[Sin94]{b:Sinai94}
{Ya}.~G. Sina{\u\i}.
\newblock {\em Topics in ergodic theory}, volume~44 of {\em Princeton
  Mathematical Series}.
\newblock Princeton University Press, Princeton, NJ, 1994.

\bibitem[Smo94]{b:Sm94}
Joel Smoller.
\newblock {\em Shock Waves and Reaction-Diffusion Equations}.
\newblock Springer-Verlag, 2nd edition, 1994.

\bibitem[Tem88]{b:Te88}
Roger Temam.
\newblock {\em Infinite Dimensional Dynamical Systems in Mechanics and
  Physics}.
\newblock Springer-Verlag, New York, 1988.

\bibitem[Ver97]{b:Veretennikov97}
A.~Yu. Veretennikov.
\newblock On polynomial mixing bounds for stochastic differential equations.
\newblock {\em Stochastic Process. Appl.}, 70(1):115--127, 1997.

\end{thebibliography}
\end{document}